\documentclass[journal]{IEEEtran}

\usepackage{cite}
\usepackage{graphicx}
\usepackage{mathtools}
\usepackage{fixltx2e}
\usepackage{algorithm}
\usepackage{url}
\usepackage[noend]{algpseudocode}
\usepackage[caption=false,font=footnotesize]{subfig}
\usepackage[absolute]{textpos}
\usepackage[hyperindex,hidelinks]{hyperref}
\usepackage[hyphenbreaks]{breakurl}

\newcommand*\Let[2]{\State #1 $\gets$ #2}
\algrenewcommand\alglinenumber[1]{
    {\sf\footnotesize\addfontfeatures{Colour=888888,Numbers=Monospaced}#1}}
\algrenewcommand\algorithmicrequire{\textbf{Input:}}
\algrenewcommand\algorithmicensure{\textbf{Output:}}
\MakeRobust{\Call}

\TPMargin{3pt}
\newcommand{\copyrightstatement}{ 
\begin{textblock}{0.84}(0.08,0.935) 
	\footnotesize
	\noindent
	\copyright~2016 IEEE. Personal use of this material is permitted. Permission from IEEE must be obtained for all other uses, in any current or future media, including reprinting/republishing this material for advertising or promotional purposes, creating new collective works, for resale or redistribution to servers or lists, or reuse of any copyrighted component of this work in other works. DOI: \href{http://dx.doi.org/10.1109/TSMC.2016.2582745}{10.1109/TSMC.2016.2582745}
\end{textblock}
}

\newcommand{\fref}[1]{Fig.~\ref{#1}}
\newcommand{\sref}[1]{Section~\ref{#1}}
\newcommand{\tref}[1]{Table~\ref{#1}}
\newcommand{\aref}[1]{Alg.~\ref{#1}}
\newcommand{\eref}[1]{(\ref{#1})}
\newcommand{\bigO}[1]{\ensuremath{\mathcal{O}\bigl(#1\bigr)}}

\newcommand{\power}{\ensuremath{P}}
\newcommand{\thrust}{\ensuremath{T}}
\newcommand{\fluidDensity}{\ensuremath{\rho}}
\newcommand{\discArea}{\ensuremath{\varsigma}}
\newcommand{\frameWeight}{\ensuremath{W}}
\newcommand{\batteryAndPayloadWeight}{\ensuremath{m}}
\newcommand{\gravity}{\ensuremath{g}}
\newcommand{\numberOfRotors}{\ensuremath{n}}
\newcommand{\helicopterPower}{\ensuremath{P^*}}
\newcommand{\singleRotorPower}{\ensuremath{\power'}}
\newcommand{\singleRotorMass}{\ensuremath{\batteryAndPayloadWeight'}}
\newcommand{\powerSlope}{\ensuremath{\alpha}}
\newcommand{\powerIntercept}{\ensuremath{\beta}}

\newcommand{\budget}{\ensuremath{B}}
\newcommand{\deliveryTimeLimit}{\ensuremath{T}}
\newcommand{\locationSet}{\ensuremath{\mathcal{N}}}
\newcommand{\locationSetWithoutDepot}{\ensuremath{\locationSet_{0}}}
\newcommand{\locationSetSum}[2]{\ensuremath{\sum_{\substack {#2 \in \locationSet \\ #1 \neq #2}}}}
\newcommand{\demand}[1]{\ensuremath{D_{#1}}}
\newcommand{\timeAtLocation}{\ensuremath{\tau}}
\newcommand{\energyDensity}{\ensuremath{\xi}}
\newcommand{\energyCostDrone}{\ensuremath{\epsilon}}
\newcommand{\droneVelocity}{\ensuremath{v}}
\newcommand{\edgeVariable}[2]{\ensuremath{x_{#1 #2}}}
\newcommand{\reuseDecision}[2]{\ensuremath{\sigma_{#1 #2}}}
\newcommand{\numberOfDrones}{\ensuremath{M}}
\newcommand{\payloadWeight}{\ensuremath{y}}
\newcommand{\payloadWeightBetweenLocations}[2]{\ensuremath{\payloadWeight_{#1 #2}}}
\newcommand{\batteryWeightBetweenLocations}[2]{\ensuremath{q_{#1 #2}}}
\newcommand{\batteryWeightAtLocation}[1]{\ensuremath{\zeta_{#1}}}
\newcommand{\largeConstant}{\ensuremath{K}}
\newcommand{\droneCarryingCapacity}{\ensuremath{Q}}
\newcommand{\locationVisitTime}[1]{\ensuremath{t_{#1}}}
\newcommand{\distance}[2]{\ensuremath{d_{#1 #2}}}
\newcommand{\depotVisitTime}[1]{\ensuremath{a_{#1}}}
\newcommand{\overallDeliveryTime}{\ensuremath{l}}
\newcommand{\fuelConsumedToReachLocation}[1]{\ensuremath{f_{#1}}}
\newcommand{\fuelConsumedToReachDepot}[1]{\ensuremath{z_{#1}}}
\newcommand{\powerFunction}[1]{\ensuremath{p(#1)}}
\newcommand{\vehicleWeightBetweenLocations}[2]{\ensuremath{m_{#1 #2}}}
\newcommand{\singleDroneCost}{\ensuremath{F}}
\newcommand{\totalCost}{\ensuremath{c}}
\newcommand{\forAllDepotlessLocations}[1]{\ensuremath{\forall #1 \in \locationSetWithoutDepot}}
\newcommand{\forAllLocations}[1]{\ensuremath{\forall #1 \in \locationSet}}
\newcommand{\forAllTwoDimensional}[4]{\ensuremath{\forall(#1,#2) \in #3 \times #4}}
\newcommand{\edgeVariableUpperBound}[2]{\ensuremath{\largeConstant \left(1 - \edgeVariable{#1}{#2}\right) }}
\newcommand{\batteryTypeSet}{\ensuremath{\mathcal{B}}}
\newcommand{\batteryTypeSum}[1]{\ensuremath{\sum_{#1 \in \batteryTypeSet}}}
\newcommand{\batteryTypeEnergy}[1]{\ensuremath{E_{#1}}}
\newcommand{\batteryTypeWeight}[1]{\ensuremath{w_{#1}}}
\newcommand{\batteryTypeCost}[1]{\ensuremath{C_{#1}}}
\newcommand{\batteryTypeAssignmentForLocation}[2]{\ensuremath{b_{#1 #2}}}
\newcommand{\batteryTypeAssignment}[1]{\ensuremath{b_{#1}}}

\newcommand{\currentSolution}{\ensuremath{\mathbf{s}}}
\newcommand{\routeVector}[1]{\ensuremath{\mathbf{r}_{#1}}}
\newcommand{\routeVectorElement}[2]{\ensuremath{r_{#1#2}}}
\newcommand{\numberOfRoutes}{\ensuremath{R}}
\newcommand{\minimizeCost}{\ensuremath{\Phi}}
\newcommand{\energyCost}{\ensuremath{\lambda}}
\newcommand{\costOfDrones}{\ensuremath{\gamma}}
\newcommand{\timeCounter}{\ensuremath{t}}
\newcommand{\travelTimeBetweenLocations}[2]{\ensuremath{\phi_{#1 #2}}}
\newcommand{\weightTimeProductCounter}{\ensuremath{\omega}}
\newcommand{\currentSolutionElement}[1]{\ensuremath{s_{#1}}}
\newcommand{\batteryEnergy}{\ensuremath{E}}
\newcommand{\firstHelperSum}{\ensuremath{Y}}
\newcommand{\secondHelperSum}{\ensuremath{Z}}
\newcommand{\routeTimingVector}{\ensuremath{\mathbf{u}}}
\newcommand{\routeTimingVectorElement}[1]{\ensuremath{u_{#1}}}
\newcommand{\deliveryTimeVectorElement}[1]{\ensuremath{h_{#1}}}
\newcommand{\arrivalTimeVectorElement}[1]{\ensuremath{a_{#1}}}
\newcommand{\currentNumberOfDrones}{\ensuremath{n}}
\newcommand{\droneTimingVector}{\ensuremath{\boldsymbol{\kappa}}}
\newcommand{\droneTimingVectorElement}[1]{\ensuremath{\kappa_{#1}}}
\newcommand{\droneDeliveryTimeVectorElement}[1]{\ensuremath{\deliveryTimeVectorElement{#1}'}}
\newcommand{\droneArrivalTimeVectorElement}[1]{\ensuremath{\arrivalTimeVectorElement{#1}'}}
\newcommand{\lowerBound}{\ensuremath{\varphi}}
\newcommand{\upperBound}{\ensuremath{\psi}}
\newcommand{\currentTemperature}{\ensuremath{\Upsilon}}
\newcommand{\initialTemperature}{\ensuremath{\currentTemperature_0}}
\newcommand{\finalTemperature}{\ensuremath{\currentTemperature'}}
\newcommand{\temperatureAdjustment}{\ensuremath{\mu}}
\newcommand{\temperatureRounds}{\ensuremath{\Lambda}}
\newcommand{\randomExchangeRule}{\ensuremath{R}}
\newcommand{\randomNumber}{\ensuremath{X}}

\newcommand{\areaSize}{\ensuremath{A}}
\newcommand{\fixedBatteryWeight}{\ensuremath{b}}

\newcommand{\singleFigureScaleValue}{0.25\textwidth}
\newcommand{\doubleFigureScaleValue}{0.2\textwidth}
\newcommand{\tripleFigureScaleValue}{0.6}

\begin{document}

\title{Vehicle Routing Problems for Drone Delivery}

\author{Kevin~Dorling,~\IEEEmembership{Student~Member,~IEEE},~Jordan~Heinrichs,~Geoffrey~G.~Messier,~\IEEEmembership{Member,~IEEE}, Sebastian~Magierowski,~\IEEEmembership{Member,~IEEE}
\thanks{This work was supported by NSERC under Grant PGSD3-443041-2013.}%
\thanks{K. Dorling, J. Heinrichs, and G. G. Messier are with the Department of Electrical and Computer Engineering, University of Calgary, Calgary, AB T2N 1N4, Canada (e-mail: kudorlin@ucalgary.ca, jheinric@ucalgary.ca, gmessier@ucalgary.ca).}
\thanks{S. Magierowski is with the Department of Electrical Engineering and Computer Science, York University, Toronto, ON M3J 1P3, Canada (e-mail:magiero@cse.yorku.ca).}}

\copyrightstatement

\maketitle

\begin{abstract}
Unmanned aerial vehicles, or drones, have the potential to significantly reduce the cost and time of making last-mile deliveries and responding to emergencies. Despite this potential, little work has gone into developing vehicle routing problems (VRPs) specifically for drone delivery scenarios. Existing VRPs are insufficient for planning drone deliveries: either multiple trips to the depot are not permitted, leading to solutions with excess drones, or the effect of battery and payload weight on energy consumption is not considered, leading to costly or infeasible routes. We propose two multi-trip VRPs for drone delivery that address both issues. One minimizes costs subject to a delivery time limit, while the other minimizes the overall delivery time subject to a budget constraint. We mathematically derive and experimentally validate an energy consumption model for multirotor drones, demonstrating that energy consumption varies approximately linearly with payload and battery weight. We use this approximation to derive mixed integer linear programs for our VRPs. We propose a cost function that considers our energy consumption model and drone reuse, and apply it in a simulated annealing (SA) heuristic for finding sub-optimal solutions to practical scenarios. To assist drone delivery practitioners with balancing cost and delivery time, the SA heuristic is used to show that the minimum cost has an inverse exponential relationship with the delivery time limit, and the minimum overall delivery time has an inverse exponential relationship with the budget. Numerical results confirm the importance of reusing drones and optimizing battery size in drone delivery VRPs.
\end{abstract}


\section{Introduction}

Unmanned aerial vehicles (UAVs), or \emph{drones}, have the potential to significantly reduce the cost and time required to deliver packages. In general, drones are less expensive to maintain than traditional delivery vehicles such as trucks, and can lower labor costs by performing tasks autonomously. Delivering with drones may be faster than delivering with traditional delivery vehicles, as drones are not limited by established infrastructure such as roads, and generally face less complex obstacle avoidance scenarios. The application of drones in emergency situations has the potential to save lives: drones can transport much-needed water, food, and medical supplies over hazardous terrain during a crisis, and can deploy wireless sensors to provide immediate updates on the event. This makes drones suitable not only for parcel delivery, but also for emergency response in the event of forest fires, oil spills, and earthquakes.

A number of large organizations have shown interest in drone delivery. In 2013, Amazon announced \emph{Prime Air} \cite{primeair}, a service that utilizes multirotor drones to deliver packages from Amazon to customers. German logistics company Deutsche Post DHL also started its \emph{Parcelcopter} project in 2013; the Parcelcopter \cite{parcelcopter} has transported medicine to the island of Juist in the North Sea. Google revealed \emph{Project Wing} \cite{projectwing} in 2014 to produce drones that can deliver larger items than Prime Air and Parcelcopter. In 2014, the United Arab Emirates announced their plan to use drones to distribute official government documents such as permits and ID cards \cite{uaedrone}. A startup called Matternet has partnered with Swiss Post to test a lightweight package delivery quadcopter \cite{matternet}.

Developments in numerous technologies have made it feasible for these organizations to perform drone deliveries. Carbon fiber manufacturing costs have decreased from \$25/kg to \$10/kg over the last 20 years \cite{Morgan2005,Reiter2013}, enabling the development of strong, lightweight airframes. Lithium polymer batteries, with their relatively high energy density \cite{Reddy2010}, have improved the flight times of UAVs compared to alternative technologies such as nickel-cadmium and nickel-metal hydride. Drones typically use GPS to determine their location, and are able to take advantage of DGPS \cite{dgps} and localization techniques \cite{karpenko2015} to improve accuracy. Obstacles can be avoided through techniques such as LIDAR \cite{Merz2013} and image processing \cite{Carloni2013}. Architectures and protocols have been developed that enable drones to form ad-hoc networks \cite{Bekmezci2013} and to wirelessly communicate with other entities \cite{Li2013}.

Even though significant effort has gone into developing drone delivery technology, challenges unique to planning drone deliveries, such as the limited flight range and carrying capacity of drones, have been neglected in comparison. The flying sidekick traveling salesman problem (FSTSP) \cite{Murray2015} provides one solution by allowing delivery trucks to also deploy drones, but relies on a combined fleet of trucks and drones to make deliveries. Green VRPs (GVRPs) \cite{Erdogan2012} may account for battery-powered vehicles, but appear to be designed with trucks or other large capacity vehicles in mind. Drones have much smaller carrying capacities than trucks, and have limited flight times, meaning that in general they can only carry a small number of packages per route. The GVRP does not allow multiple trips to the depot, so as drone capacity and the maximum travel time decrease, it would likely compensate by unnecessarily purchasing more drones even though they could have been reused. We therefore propose solving drone delivery problems with a multi-trip VRP (MTVRP) \cite{Taillard1995} that compensates for each drone's limited carrying capacity by reusing drones when possible.

The main contribution of this paper is the development of MTVRPs that consider battery and payload weight when calculating energy consumption. Other MTVRPs do not appear to account for battery and payload weight, limiting their applicability in drone delivery and emergency response scenarios. In \sref{ssec:hexacopter-energy-consumption-measurements} we experimentally show that increasing battery and payload weight can noticeably increase energy consumption, which in turn reduces flight time. Optimizing battery weight is important for drones: a drone's battery consumes a large portion of its carrying capacity, meaning that increasing its battery size to extend flight time significantly reduces the capacity available for packages. Approaches that consider the effect of payload weight on energy consumption have been shown in \cite{Xiao2012} to reduce costs in capacitated VRPs compared to approaches that do not. Balancing payload weight, battery weight, and flight time are important considerations when attempting to minimize the cost or the delivery time for drone deliveries.

In addition to developing MTVRPs for drone delivery, we are the first to derive and experimentally validate a linear energy consumption model for multirotor drones. We demonstrate that energy consumption increases at an approximately linear rate with battery and payload weight. While others have developed models for large vehicles, we are not aware of any that exist for battery-powered drones. Our experiments not only validate our energy consumption model, but also demonstrate that drones consume approximately the same power regardless of being in hover or flying at a constant speed. We will show that this considerably simplifies our MTVRPs.

Our linear energy consumption model lets us derive mixed integer linear programs (MILPs) for solving the drone delivery MTVRPs. A linear energy consumption model allows for linear energy consumption constraints, which are required to derive an MILP compatible with fast commercial solvers such as CPLEX \cite{cplex}. These problems optimize the number of drones, the routes they fly, as well as their battery weight, payload weight, and energy consumption. One VRP minimizes the cost of making the deliveries and is referred to as the \emph{minimum cost drone delivery problem (MC-DDP)}. Alternatively, in the case of an emergency response situation, a company may want to minimize overall delivery time, defined as the time it takes to complete all package deliveries; we call this the \emph{minimum time drone delivery problem (MT-DDP)}. We refer to the problems collectively as the \emph{drone delivery problems (DDPs).}

As MILPs are NP-hard, large instances (in our case, scenarios with a large number of locations) may take a significant amount of time to solve. For such cases, we derive a simulated annealing (SA) heuristic for finding sub-optimal solutions to the DDPs within a limited runtime. To help delivery network planners understand how the DDPs behave under different conditions, we perform a sensitivity analysis to show how the optimal cost and overall delivery time found by the MC-DDP and MT-DDP respectively vary with the size of the area of interest, the number of locations considered, and the restrictions placed on the budget or overall delivery time. The SA heuristic is applied to demonstrate the importance of performing multiple trips and optimizing battery weight: VRPs that do not make such considerations may find costly or infeasible solutions to the DDPs.

While our DDPs allow for multiple deliveries per route, it is not clear if delivery companies such as Amazon and DHL are interested in delivering to more than one customer per route. The Amazon Prime Air website \cite{primeair} and the DHL Parcelcopter press release \cite{parcelcopter} only show drones that can deliver single packages. Even if these companies only plan to deliver to one customer per route, our VRPs can be viewed as a generalization of this strategy. Constraints can be added to our problem to limit the number of customers per route. Furthermore, our method can be used to see if the savings from considering multiple deliveries per route make up for the costs associated with implementing the strategy.

The rest of this paper is organized as follows: in \sref{sec:related-work} we review work relevant to the DDPs. \sref{sec:energy-consumption-model} describes our energy consumption model for multirotor copters. In \sref{sec:the-ddps} we provide MILP implementations of the DDPs, while in \sref{sec:sa} we propose an SA heuristic for finding sub-optimal solutions to the problems within a constrained runtime. \sref{sec:results} presents experimental results validating the fuel consumption model, evaluates the performance of our simulated annealing heuristic, and analyzes numerical results to gain insight on the DDPs. Our conclusions are provided in \sref{sec:conclusion}.

\section{Related Work}
\label{sec:related-work}

Drone routing papers, such as \cite{Sundar2014} and \cite{Guerriero2014}, typically focus on surveillance applications. The FSTSP \cite{Murray2015} considers a delivery scenario based on combining drone delivery with a delivery truck. There do not appear to be papers that focus solely on using drones to make deliveries. This existing work on drone routing ignores factors important to drone delivery such as vehicle capacity, battery weight, changing payload weight, and reusing vehicles to reduce costs.

VRPs \cite{Cordeau2007,Eksioglu2009} have been applied to solve delivery problems that, at first glance, appear similar to the DDPs proposed in this paper. A VRP attempts to find the optimal routes for one or more vehicles to deliver commodities to a set of locations. Each location may have a unique demand representing the number or size of commodities it requires, or a time window in which the vehicle should arrive. Vehicles typically leave from one or more depots, deliver their commodities, and then return to the depots. This section will explain why existing VRPs are not adequate for the drone delivery problem.

The Green VRP (GVRP) \cite{Erdogan2012}, first described by Erdo\v{g}an et al., allows vehicles to visit refuelling stations while making deliveries to extend their range. Schneider et al. \cite{Schneider2014} created a version for battery-powered delivery vehicles that considers time windows, capacity constraints, and location demands. Hiermann et al. \cite{Hiermann2016} built on Schneider et al.'s work by allowing for a mixed fleet of various types of vehicles, perhaps with different capacities, battery sizes, and costs. They assume that charging infrastructure is in place, so vehicles can regain energy to extend their travel distance. While it could increase drone range, it is debatable whether delivery companies will deploy charging infrastructure for their drones, and unlikely that such infrastructure would be available in emergency response situations. In GVRPs, vehicles can visit the same charging station multiple times to restore energy, but appear to be able to visit the depot only once. Drones tend to have limited carrying capacities, so if the depot can only be visited once, a GVRP will probably require a large number of drones to satisfy every customer's demand. In \sref{ssec:drone-reuse} we show that reusing drones by allowing them to make multiple trips can significantly reduce costs.

A VRP that reuses vehicles is known as a multi-trip vehicle routing problem (MTVRP), first proposed by B. Fleischmann \cite{Fleischmann1990}. While a number of approaches are available for solving the MTVRP, such as a large neighborhood search algorithm \cite{Azi2014}, a hybrid genetic algorithm with a local search operator \cite{Cattaruzza2014}, a variable neighborhood search algorithm \cite{Cheikh2015}, and a branch-and-price algorithm \cite{Hernandez2016}, none appear to consider energy consumption as a function of vehicle weight. Routes found by the above approaches may be infeasible or unnecessarily costly in a drone delivery scenario, with batteries that are larger than necessary or too heavy to carry. In \sref{ssec:battery-optimization} we demonstrate that optimizing battery weight can reduce costs. Battery weight cannot be optimized without modeling energy consumption as a function of vehicle weight: the model is used to find the battery energy, and therefore weight, required to complete each route. The DDPs proposed in this paper apply our linear energy consumption model to optimize battery weight and payload weight, ensuring that the routes found are low-cost and feasible.

Our DDPs consider payload weight similar to the energy minimizing vehicle routing problem introduced by Kara et al. in \cite{Kara2007}. The authors of \cite{Kara2007} utilize a VRP that factors in the effect of a vehicle's payload on total costs, but is based solely on Newtonian physics, and is not verified against an actual vehicle. Unlike \cite{Kara2007}, the work by Xiao et al. \cite{Xiao2012} provides a linear fuel consumption model for trucks based on actual fuel consumption measurements from \cite{JapanFuelConsumption}, along with an SA algorithm for minimizing the cost of routes. The load dependant VRP \cite{Zachariadis2015} provides a local search heuristic for doing the same in pickup and delivery scenarios. While these VRPs consider the effect of payload weight on energy consumption, they appear to have been designed with vehicles that have relatively large capacities in mind, as they only let vehicles depart from the depot once. Drones have a limited carrying capacity, so if they can only leave the depot once, a large number will be required to satisfy demand in scenarios with many customers. To reduce the number and therefore cost of drones, we reuse them after they return to the depot.

The DDPs proposed in this paper combine MTVRPs with VRPs that model energy consumption as a function of vehicle weight to gain the strengths of both approaches. As mentioned earlier, MTVRPs reduce drone costs by reusing drones after they return to the depot, but can propose routes that are costly or infeasible. VRPs that model energy consumption as a function of vehicle weight, on the other hand, ensure that individual routes are feasible and low-cost, but do not reuse vehicles to lower costs. The DDPs reduce drone costs by reusing vehicles \emph{and} ensure low-cost, feasible routes by modeling energy consumption as a function of battery and payload weight.

\section{Multirotor Helicopter Energy Consumption}
\label{sec:energy-consumption-model}

The flight time of a drone is limited by its weight and the energy stored in its battery. An energy consumption model helps balance the two by providing the energy consumed by the drone as a function of its weight. When optimizing deliveries, such a model can be used to compare the energy consumed by alternative routes.

While energy consumption equations are available for single rotor helicopters \cite{Leishman2002}, to the best of our knowledge no such equations are available for multirotor helicopters, the type of drone that most delivery companies appear to favor. In this section, we derive an equation for the power consumed by a multirotor helicopter in hover as a function of its weight. We show that the power it consumes is approximately linearly proportional to the  weight of its battery and payload under practical assumptions.

We derive the power consumed by a multirotor drone in hover, but not during flight, takeoff, or landing. In flight, the power consumed by the helicopter is often reduced due to translational lift \cite{Leishman2002}, a phenomena where air flowing horizontally along the rotor generates additional lift. The average power during hover is consequently an upper bound on the average power during flight. We assume that the power consumed during takeoff and landing is, on average, approximately equivalent to the power consumed during hover. We confirm in \sref{ssec:hexacopter-energy-consumption-measurements} that these assumptions are realistic with a 3D Robotics ArduCopter Hexa-B hexacopter.

From (2.16) in \cite{Leishman2002}, we can calculate the power $\power^*$ in Watts of a single rotor helicopter in hover, with the thrust $\thrust$ in Newtons, fluid density of air $\fluidDensity$ in kg/m$^3$, and the area $\discArea$ of the spinning blade disc in m$^2$ using
\begin{equation}
\label{eq:heli-power-vs-thrust}
\helicopterPower = \frac{ \thrust^{3/2} }{\sqrt{2 \fluidDensity \discArea}},
\end{equation}
where the thrust $\thrust = (\frameWeight+\batteryAndPayloadWeight) \gravity$, given the frame weight $\frameWeight$ in kg, the battery and payload weight $\batteryAndPayloadWeight$ in kg, and gravity $\gravity$ in N.

We can use \eref{eq:heli-power-vs-thrust} to derive an equation for the power consumed by an $\numberOfRotors$-rotor copter if we assume that each rotor shares the total mass $\frameWeight+\batteryAndPayloadWeight$ of the copter equally. Each rotor carries a weight of $\singleRotorMass = \batteryAndPayloadWeight / \numberOfRotors$ for batteries and payload, and a frame weight of $\frameWeight' = \frameWeight / \numberOfRotors$, meaning that the power consumed by a single rotor is 
\begin{equation*}
\singleRotorPower = \left(\frameWeight' + \singleRotorMass \right)^{3/2} \sqrt{\frac{\gravity^3}{2 \fluidDensity \discArea}},
\end{equation*}
so the power consumed by all $\numberOfRotors$ rotors is
\begin{equation}
\label{eq:multi-power-model}
\power = \numberOfRotors \singleRotorPower = \left(\frameWeight + \batteryAndPayloadWeight \right) ^ {3/2} \sqrt{\frac{\gravity^3}{2 \fluidDensity \discArea \numberOfRotors}}.
\end{equation}
Note that the inverse relationship between $\power$ and $\numberOfRotors$ seen in \eref{eq:multi-power-model} is a result of $\numberOfRotors$ increasing the effective disc area.

The linear approximation of the $\numberOfRotors$-rotor copter power consumption equation \eref{eq:multi-power-model} can be expressed as
\begin{equation}
\label{eq:power-consumption}
\powerFunction{ \batteryAndPayloadWeight } = \powerSlope \batteryAndPayloadWeight + \powerIntercept,
\end{equation}
where $\powerSlope$ represents the power consumed per kilogram of battery and payload weight $\batteryAndPayloadWeight$, and $\powerIntercept$ is the power required to keep the hexacopter frame in the air. In \sref{sec:the-ddps} we use \eref{eq:power-consumption} to develop MILPs for optimizing the routes flown by a set of drones to minimize the time and cost required to deliver packages to a set of locations. We will see that the linear approximation is necessary for the constraints in the MILPs to be linear. Linear constraints ensure compatibility with the majority of mixed integer program solvers available.

\begin{figure}[t!]
  \centering
  \includegraphics[width=\singleFigureScaleValue]{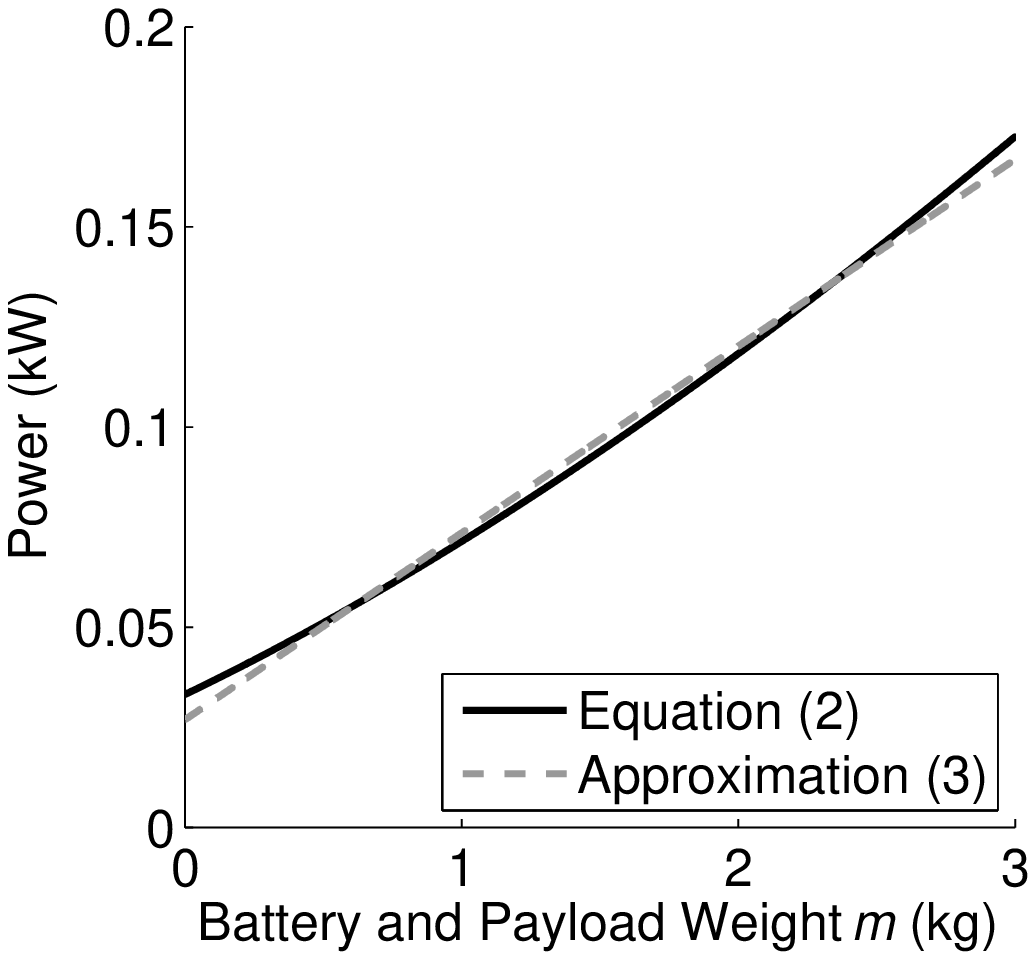}
  \caption{The linear approximation \eref{eq:power-consumption} fitted to the power consumption model \eref{eq:multi-power-model}, assuming $\numberOfRotors=6$ rotors, a fluid density of $\fluidDensity=1.204$\,kg/m$^3$, a rotor disc area of $\discArea=0.2$\,m$^2$, and a frame weight $\frameWeight=1.5$\,kg.}
  \label{fig:power-consumption}
\end{figure}

\fref{fig:power-consumption} shows that \eref{eq:power-consumption} closely fits \eref{eq:multi-power-model}, assuming a Hexa-B hexacopter such that $\numberOfRotors=6$, $\fluidDensity=1.204$\,kg/m$^3$, $\discArea=0.2$\,m$^2$, and $\frameWeight=1.5$\,kg. Applying a linear regression to \eref{eq:multi-power-model} for $\batteryAndPayloadWeight=$\,0-3\,kg in increments of 0.001\,kg results in $\powerSlope=46.7$\,W/kg and $\powerIntercept=26.9$\,W. Our linear approximation \eref{eq:power-consumption} closely fits the exact equation \eref{eq:multi-power-model}, with a mean percent error of 3.1\% and the largest difference being 6.3\,W. In \sref{ssec:hexacopter-energy-consumption-measurements} we validate this approximation experimentally with a Hexa-B hexacopter.

The linear approximation will likely hold for other multirotor helicopters. The parameters $\numberOfRotors$, $\discArea$, $\fluidDensity$, and $\frameWeight$ are presumably similar between helicopters and are all under a square root in \eref{eq:multi-power-model}, meaning slight adjustments to them will have a limited impact on $\power$. A larger variation in $\batteryAndPayloadWeight$, however, can significantly reduce the accuracy of the approximation: if $\batteryAndPayloadWeight$ varies between 0\,kg and 10\,kg, the mean percent error becomes 12.8\%, and the largest difference between \eref{eq:multi-power-model} and \eref{eq:power-consumption} is about 51\,W.

\section{The Drone Delivery Problems}
\label{sec:the-ddps}

In this section, we apply \eref{eq:power-consumption} to develop mixed integer linear programs (MILPs) for optimizing the number of drones in a fleet, as well as the routes they fly, in order to solve the drone delivery problems (DDPs). When the objective is to minimize costs, we refer to the problem as the minimum cost drone delivery problem (MC-DDP). A company that specializes in delivering parcels may prioritize cost minimization over time. In an emergency response situation, however, a company could instead prioritize time, and would consequently solve the minimum time drone delivery problem (MT-DDP). In both DDPs, constraints on the budget $\budget$ and delivery time limit $\deliveryTimeLimit$ may be imposed.

\begin{figure}[t!]
  \centering
      \includegraphics[width=\singleFigureScaleValue]{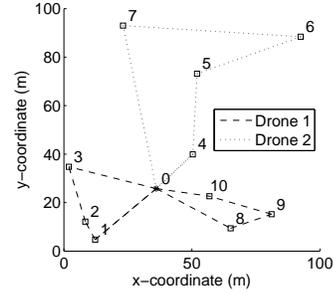}
  \caption{An example showing two drones flying three different routes.}
  \label{fig:route-example}
\end{figure}

We define a \emph{route} as a path that starts and ends at the depot. The example given in \fref{fig:route-example} depicts two drones flying three routes. Whilst Drone 2 flies a single route, visiting locations 4, 5, 6, and 7, Drone 1 flies two routes. During its first route, Drone 1 visits locations 1, 2, and 3, then returns to the depot to swap batteries and pick up more packages. After the swap, it flies its second route, visiting locations 8, 9, and 10. We assume that only one depot exists in the area. The depot acts like a charging station: each vehicle may return to the depot multiple times to collect more packages and replace its batteries. The DDPs optimize the routes flown by a set of drones in order to deliver packages to a set of locations $\locationSet$. Location 0 is the depot. Each location $i \in \locationSetWithoutDepot$, where the set $\locationSetWithoutDepot = \locationSet\ \backslash\ 0$, has a demand $\demand{i}$ which represents the weight of the package in kg that will be delivered to location $i$. Every location, except for the depot, is visited only once by a drone, and $\timeAtLocation$\,s is spent at each location to descend, deliver the package, and ascend.

The rest of this section is dedicated to explaining our assumptions and deriving the MILP forms of the MT-DDP and MC-DDP. We discuss our assumptions in \sref{ssec:assumptions}. In \sref{ssec:constraints} we list the constraints and decision variables used in the DDPs. \sref{ssec:mixed-integer-program-formulation-ddp} provides each DDP's objective function and MILP formulation.

\subsection{Assumptions}
\label{ssec:assumptions}

This section explains the assumptions we make regarding drone flight speed, battery charging, and the number of depots. We assume that drones fly between locations at a constant speed $\droneVelocity$ in m/s. In \sref{ssec:hexacopter-energy-consumption-measurements} we demonstrate that while the energy consumed by a hexacopter during flight is approximately the same as the energy consumed during hover, the energy consumed during flight can actually be slightly lower than the energy consumed in hover, likely due to translational lift. We therefore assume that the operator would fly drones as fast as they reasonably can, at a constant speed of $\droneVelocity$\,m/s between locations.

Note that one weakness of this assumption is that it ignores the impact of weather. Flying with the wind, for example, could reduce energy consumption. Cold temperatures may adversely affect battery performance until the batteries warm up. For the sake of generality, we ignore the effects of weather. Understanding how hexacopter speed and energy consumption are affected by the weather would likely require additional measurements in different weather conditions.

We assume that the demand at each location can be fully satisfied by a single drone. If the demand at a location is higher than the drone carrying capacity, the DDPs do not make multiple drones deliver optimally to that location. To overcome this issue, multiple locations can be placed on top of one another. Their combined demand can then be higher than the drone's carrying capacity. While this does not necessarily find the optimal split, it at least lets drones complete the delivery.

We do not consider recharging batteries after they have been swapped out of a drone when calculating energy costs. Instead we assume that the operator will purchase enough fully charged batteries to satisfy drone energy requirements before deliveries commence. For a delivery company, managing the charging of various battery sizes between trips could be difficult and costly compared to allowing the batteries to charge overnight. In addition, batteries have a limited number of charge cycles, so charging them between trips would likely not provide long-term savings, as they would have to be replaced more frequently. Disaster management scenarios require low response times, limiting charging time between trips.

We do not consider multiple depots. While doing so could extend the range of the drones and potentially save money, we believe considering multiple depots for drone delivery would be handled better in a separate paper that discusses optimizing the number of depots, their distribution over an area, as well as the allocation of batteries and packages among them. The survey paper \cite{MontoyaTorres2015} discusses VRPs for multiple depots.

\subsection{Decision Variables and Constraints}
\label{ssec:constraints}

This section discusses the decision variables, constants, and constraints used in the DDPs. The DDPs have identical constraints, which can be organized into categories related to limiting each drone's route, reusability, timing, energy consumption, capacity, as well as the total cost of making deliveries. Decision variables and constants are described below their corresponding constraints.

We ensure that every route is valid through
\begin{subequations}
\label{eq:flow-constraints}
\begin{align}
& \locationSetSum{i}{j} \edgeVariable{i}{j} = 1 & \forAllDepotlessLocations{i} \label{eq:visit-all-locations} \\
& \locationSetSum{i}{j} \edgeVariable{i}{j} - \locationSetSum{i}{j} \edgeVariable{j}{i} = 0 & \forAllLocations{i}. \label{eq:route-flow}
\end{align}
\end{subequations}
To create a routing map for the drones, we use the edge variables $\edgeVariable{i}{j}$, where $\edgeVariable{i}{j} = 1$ if the drone moves from location $i$ to location $j$, and $\edgeVariable{i}{j} = 0$ otherwise. Constraint \eref{eq:visit-all-locations} guarantees that every location, except for the depot, is visited exactly once by a drone, while \eref{eq:route-flow} ensures that a drone arriving at location $i$ also departs from location $i$.
 
The reusability constraints
\begin{subequations}
\label{eq:vehicle-reuse}
\begin{align}
& \sum_{j \in \locationSetWithoutDepot} \reuseDecision{i}{j} \leq \edgeVariable{i}{0} & \forAllDepotlessLocations{i} \label{eq:vehicle-reuse-1} \\
& \sum_{j \in \locationSetWithoutDepot} \reuseDecision{j}{i} \leq \edgeVariable{0}{i} & \forAllDepotlessLocations{i} \label{eq:vehicle-reuse-2} \\
& \sum_{i \in \locationSetWithoutDepot} \edgeVariable{0}{i} - \sum_{\mathclap{ \substack{(i,j) \in \locationSetWithoutDepot \times \locationSetWithoutDepot \\ i \neq j}}} \reuseDecision{i}{j} \leq \numberOfDrones & \label{eq:limit-number-of-drones}
\end{align}
\end{subequations}
determine whether or not a drone can be reused after returning to the depot. The reuse decision variable is $\reuseDecision{i}{j}$, where $\reuseDecision{i}{j}=1$ if the drone leaves location $i$ for the depot, gains a fresh battery and set of packages, then flies to location $j$ to begin a new route; otherwise $\reuseDecision{i}{j}=0$. Constraint \eref{eq:vehicle-reuse-1} implies that if a drone returns to the depot from location $i$, it is available for use again to fly to another location. Constraint \eref{eq:vehicle-reuse-2} ensures that if a reused drone leaves from the depot to location $i$, it arrived previously from another location. The number of drones that can be purchased, and can therefore fly simultaneously, is limited to $\numberOfDrones$ by \eref{eq:limit-number-of-drones}.

We apply the demand constraints
\begin{subequations}
\label{eq:demand-constraints}
\begin{align}
& \locationSetSum{i}{j} \payloadWeightBetweenLocations{j}{i} - \locationSetSum{i}{j} \payloadWeightBetweenLocations{i}{j} = \demand{i} & \forAllDepotlessLocations{i} \label{eq:demand-flow} \\
& \payloadWeightBetweenLocations{i}{j} \leq \largeConstant \edgeVariable{i}{j} & \forAllTwoDimensional{i}{j}{\locationSet}{\locationSet} \label{eq:link-demand-to-path}\\
& & i \neq j \nonumber
\end{align}
\end{subequations}
to ensure that each location receives what it demands. The payload weight between locations $i$ and $j$ is represented by the decision variable $\payloadWeightBetweenLocations{i}{j}$ in kg. The constant $\largeConstant$ is a large value representing an upper bound for constraints. Constraint \eref{eq:demand-flow} makes sure that the payload weight when leaving location $i$ is $\demand{i}$\,kg less than upon arrival. Constraint \eref{eq:link-demand-to-path} sets the payload weight of each edge without a vehicle to 0\,kg.

We calculate and enforce timing through
\begin{subequations}
\label{eq:timing-constraints}
\begin{align}
& \locationVisitTime{i} - \locationVisitTime{j} + \timeAtLocation + \distance{i}{j}/\droneVelocity & \forAllTwoDimensional{i}{j}{\locationSet}{\locationSetWithoutDepot} \label{eq:count-time} \\
& \qquad \leq \edgeVariableUpperBound{i}{j} & i \neq j \nonumber \\
& \locationVisitTime{i} - \depotVisitTime{i} + \timeAtLocation + \distance{i}{0}/\droneVelocity & \forAllDepotlessLocations{i} \label{eq:count-route-time} \\
& \qquad \leq \edgeVariableUpperBound{i}{0} & \nonumber \\
& \depotVisitTime{i} - \locationVisitTime{j} + \timeAtLocation + \distance{0}{j}/\droneVelocity & \forAllTwoDimensional{i}{j}{\locationSetWithoutDepot}{\locationSetWithoutDepot} \label{eq:connect-route-times} \\
& \qquad \leq \largeConstant \left(1 - \reuseDecision{i}{j}\right) & i \neq j \nonumber \\
& \locationVisitTime{i} \leq \overallDeliveryTime & \forAllDepotlessLocations{i} \label{eq:total-time} \\
& \overallDeliveryTime \leq \deliveryTimeLimit. \label{eq:constrain-time}
\end{align}
\end{subequations}
A location $i \in \locationSetWithoutDepot$ is visited by a drone at time $\locationVisitTime{i}$ in seconds. The time in seconds that a drone returns to the depot directly after leaving location $i$ is $\depotVisitTime{i}$. Note that $\depotVisitTime{i}=0$ if $\edgeVariable{i}{0}=0$, and $\depotVisitTime{i}>0$ otherwise. The overall delivery time $\overallDeliveryTime$ is the time in seconds when all drones have completed their deliveries. Constant values related to timing include the speed $\droneVelocity$ in m/s of the drones in the air, the distance $\distance{i}{j}$ in m between locations $i$ and $j$, as well as the time $\timeAtLocation$ in seconds spent at each location descending, delivering a package, and ascending. Drones must complete their deliveries by the delivery time limit $\deliveryTimeLimit$ in seconds. Constraint \eref{eq:count-time} keeps track of the time $\locationVisitTime{i}$ that each location $i$ is visited by a drone. Similarly, \eref{eq:count-route-time} keeps track of the time $\depotVisitTime{i}$ that a drone arrives at the depot from location $i$. Constraint \eref{eq:connect-route-times} ensures that times are correct for drones that are reused after returning to the depot. The overall delivery time $\overallDeliveryTime$ is set by \eref{eq:total-time}, and the delivery time limit $\deliveryTimeLimit$ is guaranteed by \eref{eq:constrain-time}.

Carrying capacity is restricted through
\begin{subequations}
\label{eq:capacity-constraints}
\begin{align}
& \batteryWeightBetweenLocations{i}{j} + \payloadWeightBetweenLocations{i}{j} \leq \droneCarryingCapacity \edgeVariable{i}{j} & \forAllTwoDimensional{i}{j}{\locationSet}{\locationSet} \label{eq:constrain-capacity} \\
& & i \neq j \nonumber \\
& \fuelConsumedToReachDepot{i}/\energyDensity - \batteryWeightAtLocation{i} \leq \largeConstant \left( 1 - \edgeVariable{i}{0} \right) & \forAllDepotlessLocations{i} \label{eq:set-location-energy-to-route-energy} \\
& \batteryWeightAtLocation{i} - \batteryWeightAtLocation{j} \leq \largeConstant \left( 1 - \edgeVariable{j}{i} \right) & \forAllTwoDimensional{i}{j}{\locationSetWithoutDepot}{\locationSetWithoutDepot} \label{eq:battery-capacity-at-location} \\
& & i \neq j \nonumber \\
& \batteryWeightBetweenLocations{i}{j} \geq \batteryWeightAtLocation{j} - \largeConstant \left( 1 - \edgeVariable{i}{j} \right) & \forAllTwoDimensional{i}{j}{\locationSet}{\locationSetWithoutDepot} \label{eq:battery-weight-on-path} \\
& & i \neq j \nonumber \\
& \batteryWeightBetweenLocations{i}{0} \geq \batteryWeightAtLocation{i} - \largeConstant \left( 1 - \edgeVariable{i}{0} \right) & \forAllDepotlessLocations{i}. \label{eq:battery-weight-on-path-2}
\end{align} 
\end{subequations}
The battery weight between locations $i$ and $j$ is represented by the decision variable $\batteryWeightBetweenLocations{i}{j}$ in kg. To assist with optimizing the battery weight, the decision variable $\batteryWeightAtLocation{i}$ tracks the battery weight in kg at location $i$. The energy consumed from a drone's battery by the time it arrives at the depot directly after leaving location $i$ is $\fuelConsumedToReachDepot{i}$\,kJ. Note that $\fuelConsumedToReachDepot{i}=0$ if $\edgeVariable{i}{0}=0$, while $\fuelConsumedToReachDepot{i}>0$ otherwise. The constant $\energyDensity$ is the energy density of the battery in kJ/kg. The capacity of the drone between locations $i$ and $j$ is restricted to $\droneCarryingCapacity$\,kg by \eref{eq:constrain-capacity}. The weight $\batteryWeightAtLocation{i}$ of the battery at each location $i$ is set by \eref{eq:set-location-energy-to-route-energy} and \eref{eq:battery-capacity-at-location}: constraint \eref{eq:set-location-energy-to-route-energy} finds $\batteryWeightAtLocation{i}$ for the locations visited just before the depot, while \eref{eq:battery-capacity-at-location} sets $\batteryWeightAtLocation{i}=\batteryWeightAtLocation{j}$ if the drone flies directly from location $i$ to location $j$. The weight $\batteryWeightBetweenLocations{i}{j}$ of the battery between locations $i$ and $j$ is required by \eref{eq:constrain-capacity} and is found through constraints \eref{eq:battery-weight-on-path} and \eref{eq:battery-weight-on-path-2}. Constraint \eref{eq:battery-weight-on-path} sets $\batteryWeightBetweenLocations{i}{j} \geq \batteryWeightAtLocation{j}$ if the drone flies between locations $i$ and $j$, while \eref{eq:battery-weight-on-path-2} sets $\batteryWeightBetweenLocations{i}{0} \geq \batteryWeightAtLocation{i}$ if the drone flies from location $i$ to the depot.

Energy restrictions are enforced by
\begin{subequations}
\label{eq:energy-constraints}
\begin{align}
& \fuelConsumedToReachLocation{i} - \fuelConsumedToReachLocation{j} + \powerFunction{ \vehicleWeightBetweenLocations{i}{j} }\left( \distance{i}{j} / \droneVelocity + \timeAtLocation \right) &  \forAllTwoDimensional{i}{j}{\locationSet}{\locationSetWithoutDepot} \label{eq:count-energy} \\
& \qquad \leq \edgeVariableUpperBound{i}{j} & i \neq j \nonumber \\
& \fuelConsumedToReachLocation{i} - \fuelConsumedToReachDepot{i} + \powerFunction{ \vehicleWeightBetweenLocations{i}{0} }\left( \distance{i}{0} / \droneVelocity + \timeAtLocation \right) & \forAllDepotlessLocations{i} \label{eq:route-energy} \\
& \qquad \leq \edgeVariableUpperBound{i}{0} & \nonumber \\
& \fuelConsumedToReachDepot{i} \leq \largeConstant \edgeVariable{i}{0} & \forAllDepotlessLocations{i}. \label{eq:route-energy-2}
\end{align}
\end{subequations}
The decision variable $\fuelConsumedToReachLocation{i}$ represents the energy in kJ consumed from a drone's current battery upon reaching location $i \in \locationSetWithoutDepot$. The power $\powerFunction{ \batteryAndPayloadWeight }$ in kW consumed by a drone with a battery and payload weight of $\batteryAndPayloadWeight$\,kg is the linear approximation \eref{eq:power-consumption}. The weight in kg of the drone's battery and payload between locations $i$ and $j$ is $\vehicleWeightBetweenLocations{i}{j} = \batteryWeightBetweenLocations{i}{j} + \payloadWeightBetweenLocations{i}{j}.$ Constraint \eref{eq:count-energy} forces $\fuelConsumedToReachLocation{i}$ to equal the total energy consumed along the route up to location $i$. Constraint \eref{eq:route-energy} makes $\fuelConsumedToReachDepot{i}$ equal to the energy consumed flying the entire route that ends at location $i$. To ensure that $\fuelConsumedToReachDepot{i}=0$ if the drone does not fly from location $i$ to the depot, and that $\fuelConsumedToReachDepot{i}>0$ otherwise, we include \eref{eq:route-energy-2}. Note that constraints \eref{eq:count-energy} and \eref{eq:route-energy} are linear because the power $\powerFunction{ \vehicleWeightBetweenLocations{i}{j} }$ is the linear approximation \eref{eq:power-consumption}.

Costs are kept in line with the budget $\budget$ through
\begin{subequations}
\label{eq:cost-constraints}
\begin{align}
& \totalCost = \singleDroneCost \sum_{i \in \locationSetWithoutDepot} \edgeVariable{0}{i} - \singleDroneCost \sum_{\mathclap{ \substack{(i,j) \in \locationSetWithoutDepot \times \locationSetWithoutDepot \\ i \neq j}}} \reuseDecision{i}{j} + \energyCostDrone \sum_{i \in \locationSetWithoutDepot} \fuelConsumedToReachDepot{i} \label{eq:total-costs} \\
& \totalCost \leq \budget. \label{eq:constrain-costs}
\end{align}
\end{subequations}
The cost of a drone is $\singleDroneCost$ financial units, while the budget is limited to $\budget$ financial units. The constant $\energyCostDrone$ is the cost in financial units of a kJ of energy. Constraint \eref{eq:total-costs} calculates the total cost $\totalCost$ of performing the deliveries. The leftmost term of \eref{eq:total-costs} represents the cost of drones assuming that each route requires a new drone, while the second term represents the savings provided by reusing drones; together they equal the total cost of drones. The rightmost term is the cost of energy. The total cost is restricted to the budget by \eref{eq:constrain-costs}.

The constraints in this section assume each drone's battery is sized to provide exactly enough energy for the upcoming route. The constraints can, however, be adjusted to find the optimum combination of discrete battery sizes. Assume a set of battery types $\batteryTypeSet$ exists, where each battery type $j\in\batteryTypeSet$ has energy $\batteryTypeEnergy{j}$ in kJ, a cost $\batteryTypeCost{j}$ in financial units, a weight $\batteryTypeWeight{j}$ in kg, and a decision variable $\batteryTypeAssignmentForLocation{j}{i}$ that is 1 if battery type $j$ is in the drone at location $i$ with $\fuelConsumedToReachDepot{i}\ge0$, and 0 otherwise. In \eref{eq:set-location-energy-to-route-energy} the continuous battery weight $\fuelConsumedToReachDepot{i}/\energyDensity$ at location $i$ can be replaced with the weight of the chosen batteries $\batteryTypeSum{j}\batteryTypeWeight{j}\batteryTypeAssignmentForLocation{j}{i}$. In \eref{eq:total-costs} the continuous total cost $\energyCostDrone \sum_{i \in \locationSetWithoutDepot} \fuelConsumedToReachDepot{i}$ of batteries can be replaced with the total cost of the chosen batteries $\sum_{i \in \locationSetWithoutDepot}\batteryTypeSum{j}\batteryTypeCost{j}\batteryTypeAssignmentForLocation{j}{i}$. To ensure that the chosen batteries' energy is adequate, constraint $\batteryTypeSum{j}\batteryTypeEnergy{j}\batteryTypeAssignmentForLocation{j}{i} \ge \fuelConsumedToReachDepot{i}, \forAllDepotlessLocations{i}$ can be added.

\subsection{Mixed Integer Linear Program Formulation}
\label{ssec:mixed-integer-program-formulation-ddp}

This section provides mixed integer linear program (MILP) formulations of the DDPs. As both problems have identical decision variables and constraints, only their objective functions differ. By representing the problems in MILP form, we show that while the problems are non-convex, they are compatible with commercial MIP solvers such as CPLEX \cite{cplex}.

The MT-DDP can be expressed as
\begin{align}
\min \ & \overallDeliveryTime \label{eq:mt-ddp-objective} \\
\text{s.t.} \ & \eref{eq:visit-all-locations} \ldots \eref{eq:constrain-costs}. \nonumber
\end{align}
The objective function \eref{eq:mt-ddp-objective} minimizes the overall delivery time $\overallDeliveryTime$. According to \eref{eq:total-time}, the overall delivery time is equal to the time of the last package delivery.

The MC-DDP can be expressed as
\begin{align}
\min \ & \totalCost \label{eq:mc-ddp-objective} \\
\text{s.t.} \ & \eref{eq:visit-all-locations} \ldots \eref{eq:constrain-costs}. \nonumber
\end{align}
The objective function \eref{eq:mc-ddp-objective} minimizes the total cost $\totalCost$, and therefore, according to \eref{eq:total-costs}, the cost of drones and energy.

Both DDPs are non-convex as they contain binary decision variables. However, the objective functions \eref{eq:mt-ddp-objective} and \eref{eq:mc-ddp-objective}, and the constraints \eref{eq:visit-all-locations}-\eref{eq:constrain-costs}, are all linear functions of the decision variables. This means the DDPs are mixed integer linear programs (MILPs) compatible with most MIP solvers.

While MIP solvers find optimal solutions, the runtime can be large, even for problems with a small number of locations. In \sref{sec:sa} we propose a simulated annealing approach similar to the one in \cite{Xiao2012} for finding sub-optimal solutions to the problem under limited runtimes.

\section{Simulated Annealing Implementation}
\label{sec:sa}

Even for scenarios with a small number of locations, the time required to solve the MILP implementation of the DDPs to optimality may be prohibitively long. Under time constraints, a sub-optimal solution may be preferable. In this section, we create a function for determining the cost and overall delivery time of a DDP solution. This cost function takes into account the cost and weight of batteries, the payload weight of each drone, and the fact that each drone can perform multiple trips. To solve practical problems with large numbers of locations, we apply it to a simulated annealing (SA) \cite{Kirkpatrick1983} heuristic similar to the algorithms in \cite{Xiao2012} and \cite{Kaku2003}.

We take a \emph{string-based} approach, meaning that a DDP solution is represented as a one-dimensional vector of whole numbers
$
\currentSolution = [0 \ \routeVector{1} \ 0 \ \routeVector{2} \ 0 \ \ldots \ \routeVector{\numberOfRoutes} \ 0],
$
where, as described in \cite{Kaku2003}, the number 0 represents the depot, $\routeVector{n}$ represents a route flown by a drone, and $\numberOfRoutes$ is the number of routes. The variable $\routeVector{n}$ is a vector whose $k$th element $\routeVectorElement{n}{k} \in \locationSet$, where $\locationSet$ is the set of locations. The routes in \fref{fig:route-example} could, for example, be represented by the string $\currentSolution = [0 \ 1 \ 2 \ 3 \ 0 \ 0 \ 0 \ 8 \ 9 \ 10 \ 0 \ 4 \ 5 \ 6 \ 7 \ 0]$. When two neighboring zeros occur, there is an empty route; a total of three routes exist in the above $\currentSolution$, even though five could have been flown. Section 2.3 of \cite{Kaku2003} discusses the string-based model and its advantages in greater detail.

The rest of this section is dedicated to explaining our SA implementation. To decide whether one solution is an improvement over another, we compare their costs. The method we use to determine the cost of a solution is presented in \sref{ssec:route-costs}. Finally, the simulated annealing algorithm is discussed in \sref{ssec:sa-algorithm}.

\subsection{The DDP Cost Function}
\label{ssec:route-costs}

The cost of a set of routes depends on the type of DDP being solved. In the MT-DDP the best solution has the lowest overall delivery time, while in the MC-DDP the best solution has the lowest cost. In this section, we create a function for calculating the cost and overall delivery time of a DDP solution.

We assume that violating the capacity constraint, budget constraint, or delivery time limit is costly. As in \cite{Xiao2012}, we reduce the likelihood of accepting a solution that violates a constraint by penalizing it; we do so by greatly increasing the cost and overall delivery time. The constant $\largeConstant$ can be adjusted to reflect the actual costs of violating a constraint.

\aref{alg:solution-cost} calculates the overall cost and delivery time of a solution $\currentSolution$ to the DDP. It takes the solution vector $\currentSolution$ and a boolean value $\minimizeCost$ that is 1 for the MC-DDP and 0 for the MT-DDP as inputs. The algorithm is divided into two phases: the first phase calculates the cost of energy with \aref{alg:energy-cost}, while the second finds the cost of drones and the overall delivery time with \aref{alg:drone-cost}. Later in this section we show that \aref{alg:energy-cost} has a $\bigO{|\locationSetWithoutDepot|}$ bound, and \aref{alg:drone-cost} has a $\bigO{|\locationSetWithoutDepot|\log^2\numberOfDrones}$ bound, where $|\locationSetWithoutDepot|$ is the number of locations in the area of interest, and $\numberOfDrones$ is the maximum number of drones that can be purchased. This means that \aref{alg:solution-cost} has a $\bigO{|\locationSetWithoutDepot|\log^2\numberOfDrones}$ bound as well. The rest of this section discusses Algs. \ref{alg:energy-cost}-\ref{alg:drone-cost}.

\begin{algorithm}[t!]
\caption{The \textsc{cost}$(\currentSolution,\minimizeCost)$ method.}
\label{alg:solution-cost}
\begin{algorithmic}
\Require
\Statex $\currentSolution$ - The solution vector
\Statex $\minimizeCost$ - A boolean variable that is 1 if the goal is to minimize cost, or 0 to minimize overall delivery time
\Ensure
\Statex $\totalCost$ - The total cost of $\currentSolution$
\Statex $\overallDeliveryTime$ - The time taken by $\currentSolution$ to complete deliveries
\\
\State{\emph{// Find the cost of energy $\energyCost$, the cost $\costOfDrones$ of drones, the overall}}
\State{\emph{// delivery time $\overallDeliveryTime$, and the total cost $\totalCost$.}}
\Let{$\energyCost$}{\Call{energyCost}{$\currentSolution$}}
\Let{$(\costOfDrones,\overallDeliveryTime)$}{\textsc{droneCostAndDeliveryTime}$(\currentSolution,\energyCost,\minimizeCost)$}
\Let{$\totalCost$}{$\energyCost+\costOfDrones$}
\State{\emph{// Enforce budget and delivery time limit by increasing cost}}
\State{\emph{// and overall delivery time if either constraint is violated.}}
\If{$\totalCost > \budget$}
\Let{$\totalCost$}{$\totalCost + \largeConstant \left( \totalCost - \budget \right)$}
\Let{$\overallDeliveryTime$}{$\overallDeliveryTime + \largeConstant \left( \totalCost - \budget \right)$}
\ElsIf{$\overallDeliveryTime > \deliveryTimeLimit$}
\Let{$\totalCost$}{$\totalCost + \largeConstant \left( \overallDeliveryTime - \deliveryTimeLimit \right)$}
\Let{$\overallDeliveryTime$}{$\overallDeliveryTime + \largeConstant \left( \overallDeliveryTime - \deliveryTimeLimit \right)$}
\EndIf
\State \Return ($\totalCost,\overallDeliveryTime$)
\end{algorithmic}
\end{algorithm}

\aref{alg:energy-cost} finds the cost of energy by iterating over the solution vector $\currentSolution$. It iterates from back to front, like the cost function in Fig. A2 in \cite{Xiao2012}; if it went from front to back, it would have to calculate the initial value of payload weight $\payloadWeight$ before iterating over a route, instead of setting $\payloadWeight=0$. \aref{alg:energy-cost} calls \textsc{batteryEnergy}($\timeCounter,\weightTimeProductCounter$) to find the energy
\begin{equation}
\begin{aligned}
\label{eq:battery-energy}
\batteryEnergy & = \sum_{k=a}^{b-1} \powerFunction{ \payloadWeightBetweenLocations{k}{(k+1)} + \batteryEnergy / \energyDensity} \travelTimeBetweenLocations{k}{(k+1)} \\
& = \frac{ \powerSlope \sum_{k=a}^{b-1} \payloadWeightBetweenLocations{k}{(k+1)} \travelTimeBetweenLocations{k}{(k+1)} + \powerIntercept  \sum_{k=a}^{b-1} \travelTimeBetweenLocations{k}{(k+1)} }{1 - (\powerSlope/\energyDensity) \sum_{k=a}^{b-1} \travelTimeBetweenLocations{k}{(k+1)} } \\
& = \frac{ \powerSlope \weightTimeProductCounter + \powerIntercept \timeCounter }{ 1 - (\powerSlope/\energyDensity) \timeCounter }
\end{aligned}
\end{equation}
in kJ required by a drone to complete a route containing locations $a$ to $b$. Note that $\timeCounter$ is the time taken to complete a route, $\weightTimeProductCounter$ is the sum-of-products of the travel time and payload weight between locations, $\payloadWeightBetweenLocations{i}{j}$ is the payload weight in kg between locations $i$ and $j$, $\powerFunction{\batteryAndPayloadWeight}$ is \eref{eq:power-consumption}, $\powerSlope$ is the power in kW consumed per kg of battery and payload weight, $\powerIntercept$ is the power in kW required to keep the drone in hover assuming no battery or payload weight, $\energyDensity$ is the battery energy density in kJ/kg, and $\travelTimeBetweenLocations{i}{j}$ is the travel time in seconds between locations $i$ and $j$. As \aref{alg:energy-cost} iterates through $\currentSolution$ once, its bound is $\bigO{|\locationSetWithoutDepot|}$.

\begin{algorithm}[t!]
\caption{The \textsc{energyCost}$(\currentSolution)$ method.}
\label{alg:energy-cost}
\begin{algorithmic}
\Require
\Statex $\currentSolution$ - The solution vector
\Ensure
\Statex $\energyCost$ - The cost of energy for solution $\currentSolution$
\\
\State{\emph{// Initialize the time $\timeCounter$, sum-of-products $\weightTimeProductCounter$ of time and weight,}}
\State{\emph{// payload weight $\payloadWeight$, and energy cost $\energyCost$.}}
\Let{$\timeCounter$}{0}
\Let{$\weightTimeProductCounter$}{0}
\Let{$\payloadWeight$}{$0$}
\Let{$\energyCost$}{$0$}
\State{\emph{// Iterate through the locations in $\currentSolution$ from back to front.}}
\For{$k \in \{$ \Call{length}{$\currentSolution$}-2$ \ldots 0 \}$}
\State{\emph{// Set the next location $i$ and current location $j$.}}
\Let{$i$}{$\currentSolutionElement{k+1}$}
\Let{$j$}{$\currentSolutionElement{k}$}
\If{$i \neq 0 \vee j \neq 0$}
\State{\emph{// Record the travel time $\travelTimeBetweenLocations{i}{j}$ between $i$ and $j$.}}
\Let{$\travelTimeBetweenLocations{i}{j}$}{$\timeAtLocation + \distance{i}{j} / \droneVelocity$}
\State{\emph{// Update the time, the sum-of-products of time and}}
\State{\emph{// weight, and the payload weight respectively.}}
\Let{$\timeCounter$}{$\timeCounter + \travelTimeBetweenLocations{i}{j}$}
\Let{$\weightTimeProductCounter$}{$\weightTimeProductCounter + \payloadWeight \travelTimeBetweenLocations{i}{j}$}
\Let{$\payloadWeight$}{$\payloadWeight + \demand{j}$}
\State{\emph{// Update costs when a vehicle arrives at a depot.}}
\If{$\left(0 == j\right) \wedge \left(i \neq 0 \right)$}
\Let{$\batteryEnergy$}{\Call{batteryEnergy}{$\timeCounter,\weightTimeProductCounter$}}
\State{\emph{// Add the cost of energy. Note that $\energyCostDrone$ is the cost}}
\State{\emph{// of a kJ. Greatly increase the cost if the positive}}
\State{\emph{// energy constraint is violated.}}
\If{$\batteryEnergy>0$}
\Let{$\energyCost$}{$\energyCost + \batteryEnergy \energyCostDrone $}
\Else
\Let{$\energyCost$}{$\energyCost + - \largeConstant \left(\batteryEnergy \energyCostDrone\right)$}
\EndIf
\State{\emph{// Enforce the capacity constraint by greatly}}
\State{\emph{// increasing cost if the drone carries too much.}}
\If{$\payloadWeight + \batteryEnergy / \energyDensity >\droneCarryingCapacity$}
\Let{$\energyCost$}{$\energyCost + \largeConstant \left( \payloadWeight + \batteryEnergy / \energyDensity - \droneCarryingCapacity \right)$}
\EndIf
\State{\emph{// Reset weight and time values for the next route.}}
\Let{$\timeCounter$}{0}
\Let{$\weightTimeProductCounter$}{0}
\Let{$\payloadWeight$}{$0$}
\EndIf
\EndIf
\EndFor
\Return $\energyCost$
\end{algorithmic}
\end{algorithm}

\begin{algorithm}[t!]
\caption{The \textsc{routeDeliveryAndArrivalTimes}$(\currentSolution)$ method.}
\label{alg:route-times}
\begin{algorithmic}
\Require
\Statex $\currentSolution$ - The solution vector
\Ensure
\Statex $\routeTimingVector$ - A vector containing pairs of the delivery and arrival times for each route
\\
\State{\emph{// Initialize the delivery time $\deliveryTimeVectorElement{}$ and arrival time $\arrivalTimeVectorElement{}$ to 0.}}
\Let{$\deliveryTimeVectorElement{}$}{0}
\Let{$\arrivalTimeVectorElement{}$}{0}
\State{\emph{// Initialize the route timing vector $\routeTimingVector$ to $\emptyset$.}}
\Let{$\routeTimingVector$}{$\emptyset$}
\State{\emph{// Iterate through the locations in $\currentSolution$ from front to back.}}
\For{$k \in \{1 \ldots $ \Call{length}{$\currentSolution$} $-1\}$}
\State{\emph{// Set the current location $i$ and previous location $j$.}}
\Let{$i$}{$\currentSolutionElement{k}$}
\Let{$j$}{$\currentSolutionElement{k-1}$}
\If{$i \neq 0 \vee j \neq 0$}
\State{\emph{// Update the time variables.}}
\Let{$\deliveryTimeVectorElement{}$}{$\arrivalTimeVectorElement{}$}
\Let{$\travelTimeBetweenLocations{i}{j}$}{$\timeAtLocation + \distance{i}{j} / \droneVelocity$}
\Let{$\arrivalTimeVectorElement{}$}{$\arrivalTimeVectorElement{} + \travelTimeBetweenLocations{i}{j}$}
\State{\emph{// Add the delivery and arrival time to $\routeTimingVector$.}}
\If{$i == 0 \wedge j \neq 0$}
\Let{$\routeTimingVector$}{$\routeTimingVector \cup (\deliveryTimeVectorElement{},\arrivalTimeVectorElement{})$}
\State{\emph{// Reset the delivery time value.}}
\Let{$\arrivalTimeVectorElement{}$}{0}
\EndIf
\EndIf
\EndFor
\\
\Return $\routeTimingVector$
\end{algorithmic}
\end{algorithm}

\begin{algorithm}[t!]
\caption{The \textsc{listSchedule}$(\routeTimingVector,\currentNumberOfDrones)$ method.}
\label{alg:list-schedule}
\begin{algorithmic}
\Require
\Statex $\routeTimingVector$ - A vector containing pairs of delivery and arrival times for each route
\Statex $\currentNumberOfDrones$ - The number of drones that can be assigned routes
\Ensure
\Statex $\overallDeliveryTime$ - The overall delivery time when the routes are assigned to $\currentNumberOfDrones$ drones
\\
\State \emph{// Apply the list scheduling algorithm to assign routes to}
\State \emph{// drones. Initialize the drone timing vector $\droneTimingVector$.}
\Let{$\droneTimingVector$}{$\emptyset$}
\For{$i \in \{0 \ldots \currentNumberOfDrones-1\}$}
\Let{$\droneTimingVector$}{$\droneTimingVector \cup (0,0)$}
\EndFor

\State \emph{// Go through every route in the route timing vector $\routeTimingVector$.}
\For{$i \in \{0 \ldots\ $\Call{length}{$\routeTimingVector$}$-1\}$}
\State \emph{// Get the delivery and arrival times of the drone with}
\State \emph{// the smallest arrival time. Then assign route $i$ to that}
\State \emph{//  drone.}
\Let{$(\droneDeliveryTimeVectorElement{},\droneArrivalTimeVectorElement{})$}
{\Call{popMinimumArrivalTimeElement}{$\droneTimingVector$}}
\Let{$(\deliveryTimeVectorElement{},\arrivalTimeVectorElement{})$}{$\routeTimingVectorElement{i}$}
\State \Call{insertElement}{$\droneTimingVector,(\droneArrivalTimeVectorElement{} + \deliveryTimeVectorElement{},\droneArrivalTimeVectorElement{} + \arrivalTimeVectorElement{})$}
\EndFor
\State \emph{// Return the maximum delivery time in $\droneTimingVector$.}
\\
\Return \Call{maximumDeliveryTime}{$\droneTimingVector$}

\end{algorithmic}
\end{algorithm}

\aref{alg:energy-cost} assumes that battery weight and cost can be sized to provide the exact amount of energy required. If it is preferable to combine discrete battery sizes instead, it can be shown that finding the optimal combination of batteries that minimizes cost while providing enough energy to complete the route is a minimum knapsack problem \cite{Kellerer2004}. Let $\batteryTypeSet$ represent a set of battery types, where each battery type $j\in\batteryTypeSet$ has a weight $\batteryTypeWeight{j}$ in kg, energy $\batteryTypeEnergy{j}$ in kJ, and cost $\batteryTypeCost{j}$ in financial units. Given a decision variable $\batteryTypeAssignment{j}$ that is 1 if battery type $j$ is used, and 0 otherwise, while letting $\firstHelperSum=\batteryTypeEnergy{j} - \powerSlope \sum_{k=a}^{b-1} \batteryTypeWeight{j}\travelTimeBetweenLocations{k}{(k+1)}$ and $\secondHelperSum= \powerSlope \sum_{k=a}^{b-1} \payloadWeightBetweenLocations{k}{(k+1)} \travelTimeBetweenLocations{k}{(k+1)} + \powerIntercept \sum_{k=a}^{b-1} \travelTimeBetweenLocations{k}{(k+1)}$, the minimum knapsack problem is
\begin{align}
\label{eq:min-knapsack}
\min \ & \batteryTypeSum{j}\batteryTypeCost{j}\batteryTypeAssignment{j} \\
\text{s.t.} \ & \batteryTypeSum{j}\firstHelperSum\batteryTypeAssignment{j} \ge \secondHelperSum. \nonumber
\end{align}
The values $\batteryEnergy$, $\batteryEnergy \epsilon$, and $\batteryEnergy / \xi$ in \aref{alg:energy-cost} can be replaced with $\batteryTypeSum{j}\batteryTypeEnergy{j}\batteryTypeAssignment{j}$, $\batteryTypeSum{j}\batteryTypeCost{j}\batteryTypeAssignment{j}$, and $\batteryTypeSum{j}\batteryTypeWeight{j}\batteryTypeAssignment{j}$ respectively after solving \eref{eq:min-knapsack}.

To calculate the overall delivery time, we require the \emph{delivery time} and \emph{arrival time} of each route, which can be found with \aref{alg:route-times}. The delivery time is when the last delivery along a route is made, while the arrival time is when the drone arrives at the depot after completing the route. Route timing vector $\routeTimingVector$ contains pairs of delivery and arrival times for each route, such that element $\routeTimingVectorElement{i} = (\deliveryTimeVectorElement{i},\arrivalTimeVectorElement{i})$, where $\deliveryTimeVectorElement{i}$ is the delivery time and $\arrivalTimeVectorElement{i}$ is the arrival time for route $i$. \aref{alg:route-times} implements \textsc{routeDeliveryAndArrivalTimes}$(\currentSolution)$, which iterates over the solution $\currentSolution$ to calculate $\routeTimingVector$. As \aref{alg:route-times} iterates through $\currentSolution$ once, it has a $\bigO{|\locationSetWithoutDepot|}$ bound.

To minimize the overall delivery time, we require a method of assigning $\numberOfRoutes$ routes to $\currentNumberOfDrones$ drones; this method is given in \aref{alg:list-schedule}. This problem is equivalent to the \emph{load balancing} problem in computer science, where $j$ jobs are assigned to $m$ machines to minimize the maximum load of the machines. This problem is NP-hard, so we apply the \emph{list scheduling}, or \emph{greedy-balance}, algorithm discussed in Section 11.1 of \cite{Kleinberg2011}, except jobs become routes, machines become drones, and we try to minimize the maximum delivery time of the drones. \aref{alg:list-schedule} implements \textsc{listSchedule}($\routeTimingVector$,$\currentNumberOfDrones$) to approximate the minimum overall delivery time. At most $\numberOfDrones$ drones can be purchased, and the number of elements in $\routeTimingVector$ cannot exceed the number of locations $|\locationSetWithoutDepot|$. Drone timing vector $\droneTimingVector$ contains pairs of delivery and arrival times for each route, such that element $\droneTimingVectorElement{i} = (\deliveryTimeVectorElement{i},\arrivalTimeVectorElement{i})$. Assuming that $\droneTimingVector$ is implemented as a binary heap, the bound of \aref{alg:list-schedule} is $\bigO{|\locationSetWithoutDepot|\log\numberOfDrones}$.

Alternatives to list scheduling, such as the \emph{longest processing time} (LPT), or \emph{sorted-balance}, algorithm discussed in \cite{Kleinberg2011} could be applied. LPT improves list scheduling by choosing the jobs to assign based on a specific order, such as by descending runtime, resulting in a better performance guarantee of 1.5 instead of 2. The order of the route delivery and arrival times in $\routeTimingVector$ depend on the order of the routes in $\currentSolution$. The SA algorithm randomly adjusts $\currentSolution$: in addition to changing the routes themselves, it also changes their order. As the SA algorithm is already adjusting the order of the routes, we choose not to implement LPT.

\begin{algorithm}[t!]
\caption{The \textsc{droneCostAndDeliveryTime}$(\currentSolution,\energyCost,\minimizeCost)$ method.}
\label{alg:drone-cost}
\begin{algorithmic}
\Require
\Statex $\currentSolution$ - The solution vector
\Statex $\energyCost$ - The cost of energy for solution $\currentSolution$
\Statex $\minimizeCost$ - A boolean variable that is 1 if the goal is to minimize drone cost, or 0 to minimize overall delivery time
\Ensure
\Statex $\costOfDrones$ - The cost of drones
\Statex $\overallDeliveryTime$ - The overall delivery time
\\
\State \emph{// Get the vector $\routeTimingVector$, which contains pairs of delivery and}
\State \emph{// arrival times for each route.}
\Let{$\routeTimingVector$}{\Call{RouteDeliveryAndArrivalTimes}{$\currentSolution$}}
\State \emph{// Perform a binary search to find the best number of}
\State \emph{// drones $\currentNumberOfDrones$ when minimizing cost.}
\If{$\minimizeCost$}
\State \emph{// Set the lower and upper bounds $\lowerBound$ and $\upperBound$ respectively,}
\State \emph{// where $\numberOfDrones$ is the maximum number of drones available.}
\Let{$\lowerBound$}{1}
\Let{$\upperBound$}{$\numberOfDrones$}
\While{$\lowerBound \leq \upperBound-1$}
\Let{$\currentNumberOfDrones$}{$\lfloor \lowerBound + \frac{\upperBound-\lowerBound}{2} \rfloor$}
\If{\Call{listSchedule}{$\routeTimingVector,\currentNumberOfDrones$}$\ \leq \deliveryTimeLimit$}
\Let{$\upperBound$}{$\currentNumberOfDrones$}
\Else
\Let{$\lowerBound$}{$\currentNumberOfDrones+1$}
\EndIf
\EndWhile
\Let{$\currentNumberOfDrones$}{$\lowerBound$}
\If{\Call{listSchedule}{$\routeTimingVector,\currentNumberOfDrones$}$\ > \deliveryTimeLimit$}
\Let{$\currentNumberOfDrones$}{$\upperBound$}
\EndIf
\Else
\State \emph{// Purchase as many drones as possible to minimize time,}
\State \emph{// given a budget $\budget$ and a cost $\singleDroneCost$ of a single drone.}
\Let{$\currentNumberOfDrones$}{$\lfloor \frac{\budget - \energyCost}{\singleDroneCost} \rfloor$}
\State \emph{// Make sure number of drones is positive.}
\If{$\currentNumberOfDrones < 1$}
\Let{$\currentNumberOfDrones$}{$1$}
\EndIf
\EndIf
\State \emph{// Find the cost of drones $\costOfDrones$ and the overall delivery time $\overallDeliveryTime$.}
\Let{$\costOfDrones$}{$\currentNumberOfDrones \singleDroneCost$}
\Let{$\overallDeliveryTime$}{\Call{listSchedule}{$\routeTimingVector,\currentNumberOfDrones$}}
\\
\Return $(\costOfDrones,\overallDeliveryTime)$
\end{algorithmic}
\end{algorithm}

\aref{alg:drone-cost} implements the \textsc{droneCostAndDeliveryTime}$(\currentSolution,\energyCost,\minimizeCost)$ method, where $\energyCost$ is the cost of energy, to calculate the cost of drones $\costOfDrones$ and the overall delivery time $\overallDeliveryTime$. When solving the MC-DDP, it performs a binary search to reduce the number of drones $\currentNumberOfDrones$ as much as possible, without violating the delivery time constraint. While solving the MT-DDP, on the other hand, it sets $\currentNumberOfDrones$ as high as possible without violating the budget constraint. After finding $\currentNumberOfDrones$, and in turn $\costOfDrones$, it uses \aref{alg:list-schedule} to find $\overallDeliveryTime$. As \aref{alg:drone-cost} performs a binary search to find the number of drones, calling \aref{alg:list-schedule} every iteration, it has a $\bigO{|\locationSetWithoutDepot|\log^2\numberOfDrones}$ bound.

\subsection{Simulated Annealing Algorithm}
\label{ssec:sa-algorithm}

\aref{alg:sa} provides the simulated annealing algorithm we use. It requires an initial temperature $\initialTemperature$, final temperature $\finalTemperature$, cooling factor $\temperatureAdjustment$, and the number of rounds $\temperatureRounds$ per cooling phase. It is based on the approach discussed in Section 2.5 of \cite{Kaku2003}, with slight modifications. We do not include the heating phase, and instead let $\initialTemperature$ be an input to the algorithm. We use \aref{alg:solution-cost} to compare the costs of different solutions.

\begin{algorithm}[t!]
\caption{The \textsc{simulatedAnnealing}$(\initialTemperature,\finalTemperature,\temperatureAdjustment, \temperatureRounds)$ method.}
\label{alg:sa}
\begin{algorithmic}
\Require
\Statex $\initialTemperature$ - The initial temperature
\Statex $\finalTemperature$ - The final temperature
\Statex $\temperatureAdjustment$ - The amount that temperature is adjusted each iteration
\Statex $\temperatureRounds$ - The number of adjustments to the solution attempted per iteration
\Ensure
\Statex $\currentSolution$ - A solution to the DDP
\\
\State \emph{// Generate a random solution to the DDP using the number}
\State \emph{// of locations $|\locationSetWithoutDepot|$. It can be infeasible at this point.}
\Let{$\currentSolution$}{\Call{randomSolution}{$|\locationSetWithoutDepot|$}}
\State \emph{// Set the current temperature $\currentTemperature$ to the initial temperature}
\State \emph{// $\initialTemperature$.}
\Let{$\currentTemperature$}{$\initialTemperature$}
\State \emph{// Perform annealing on the current solution $\currentSolution$.}
\While{$\currentTemperature > \finalTemperature$}
\Let{$\currentTemperature$}{$\temperatureAdjustment\currentTemperature$}
\For{$k \in \{1\ldots\temperatureRounds\}$}
\State \emph{// Pick two random indices for $\currentSolution$.}
\Let{$i$}{\Call{uniformRandomInteger}{$1$,\Call{length}{$\currentSolution$}-1}}
\Let{$j$}{\Call{uniformRandomInteger}{$1$,\Call{length}{$\currentSolution$}-1}}
\State \emph{// Choose a random exchange rule $\randomExchangeRule$, then adjust $\currentSolution$.}
\Let{$\randomExchangeRule$}{\Call{uniformRandomInteger}{$1$,$3$}}
\Let{$\currentSolution'$}{\Call{adjustWithExchangeRule}{$\currentSolution,\randomExchangeRule,i,j$}}
\State \emph{// Use Metropolis algorithm to determine if the}
\State \emph{// adjustment should become the current solution.}
\Let{$\randomNumber$}{\Call{uniformRandomFloat}{$0$,$1$}}
\If{$\exp\left(-\frac{\textsc{cost}(\currentSolution') - \textsc{cost}(\currentSolution)}{\currentTemperature}\right) \geq \randomNumber$}
\Let{$\currentSolution$}{$\currentSolution'$}
\EndIf
\EndFor
\EndWhile
\Return{\currentSolution}
\end{algorithmic}
\end{algorithm}

\begin{figure}[t!]
  \centering
  \subfloat[Swap\label{fig:exchange-rule-1}]{\includegraphics[scale=\tripleFigureScaleValue]{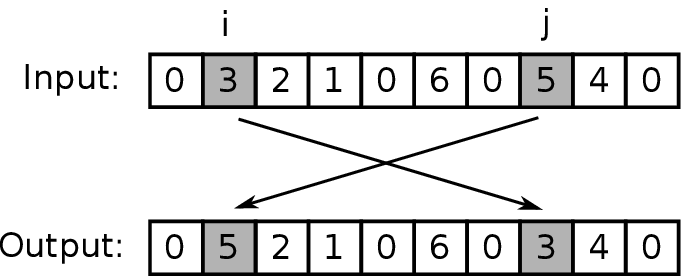}}
  \hspace{2px}
  \subfloat[Relocate\label{fig:exchange-rule-2}]{\includegraphics[scale=\tripleFigureScaleValue]{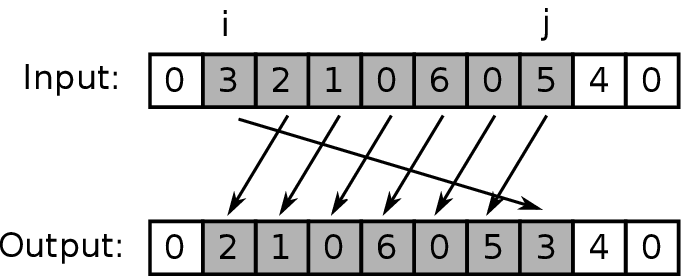}}
  \hspace{2px}
  \subfloat[2-opt\label{fig:exchange-rule-3}]{\includegraphics[scale=\tripleFigureScaleValue]{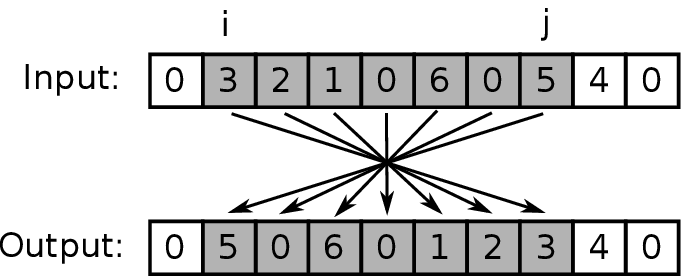}}
  \caption{The exchange rules from Section 2.4 of \cite{Kaku2003}.}
  \label{fig:exchange-rules}
\end{figure}

In \aref{alg:sa} the \emph{exchange rules} are used to find a neighboring solution $\currentSolution'$ from the current solution $\currentSolution$. The three exchange rules in Section 2.4 of \cite{Kaku2003} are applied in \aref{alg:sa}, each using two uniform randomly selected indices. The swap rule in \fref{fig:exchange-rule-1} swaps two randomly selected elements in the vector, the relocate rule in \fref{fig:exchange-rule-2} changes the position of a randomly selected element in the vector, while the 2-opt rule in \fref{fig:exchange-rule-3} reverses a randomly chosen sub-vector. Note that the vectors in \fref{fig:exchange-rules} represent the solution vector $\currentSolution$, whose elements are the locations in the area of interest. To determine if the neighboring solution should become the current one, we apply the Metropolis algorithm \cite{Metropolis1953}.

\sref{sec:results} compares the performance of the MILP implementations to the SA implementations to validate the cost function and demonstrate that the SA implementation can find near-optimal solutions to small DDP instances. One weakness of the SA approach is that it does not take advantage of characteristics unique to VRPs. For example, it does not take advantage of geographical information to reduce the likelihood of connecting distant locations, even though a route containing such locations would likely be infeasible. It also does not keep track of how often the swap, relocate, and 2-opt operations improve the solution, which could be used to increase the frequency of operations that are more likely to improve the solution. While heuristics that implement these approaches can likely find better solutions than SA given the same runtime, we do not believe it is likely that the trends observed in \sref{sec:results} will change significantly with the heuristic. The trends can be explained by adjustments made to parameters such as area size and the number of customer locations in an area, and do not appear to be a result of behavior specific to SA.

\section{Results}
\label{sec:results}

In this section, we validate the multirotor helicopter energy consumption model, and analyze the performance of the MILP and SA implementations of the DDPs. We demonstrate that the SA implementations find near-optimal solutions to small instances of the DDPs, and provide consistent results for larger instances. We show that the cost found by the MC-DDP has an inverse exponential relationship with the delivery time limit, and that the overall delivery time found by the MT-DDP has an inverse exponential relationship with the budget. We demonstrate that the majority of costs go towards drones, and that changing the budget or delivery time limit causes the SA DDPs to adjust the number of drones while leaving energy costs approximately constant. We demonstrate that reusing drones to perform multiple trips is an important strategy for reducing costs. Our results indicate that optimizing battery weight is an effective strategy for reducing the total cost and overall delivery time.

To understand how the DDPs behave in general, we average the results of a number of randomly generated instances of each scenario. Each scenario has an area size of 0.25\,km$^2$ or 1\,km$^2$. Small scenarios have 6-8 delivery locations, while large scenarios have 125 or 500 delivery locations. We generate 50 random instances for each scenario. In each instance, delivery locations are uniformly distributed throughout the area and are given uniform random demands of 0.5-2\,kg; depots are located in the middle of the area. We run the SA algorithm 20 times per instance, and calculate the minimum, mean, standard deviation, and average runtime of these 20 runs. The minimums, means, etc. in this section are averaged over all 50 instances of each scenario. Unless mentioned otherwise, we use the average minimum values in tables and graphs. We run the MILP implementation once per instance, and provide the average result and runtime over all 50 instances of a scenario. Unless mentioned otherwise, when running the SA algorithm, the initial temperature $\initialTemperature=1$, the final temperature $\finalTemperature=0.001$, the cooling factor $\temperatureAdjustment=0.99$, and the number of rounds $\temperatureRounds=1000$.

We assume that each drone costs $\singleDroneCost=$\ \$500 and has a maximum carrying capacity of $\droneCarryingCapacity=3$\,kg. Its maximum carrying capacity includes its battery and payload weight, but not its frame weight. Drones can vary in price from under \$100 for basic models to tens of thousands of dollars for industrial or military use. We assume that the operator either manufactures its own drones, or uses a hobbyist drone to make deliveries, hence the \$500 cost. A carrying capacity of 3\,kg should be possible for drones designed to deliver small packages. We assign a slope $\powerSlope=0.217$\,kW/kg and a y-intercept $\powerIntercept=0.185$\,kW to the energy consumption model, based on results obtained from the hexacopter power consumption measurements made in \sref{ssec:hexacopter-energy-consumption-measurements}. When delivering packages, we assume that the drone flies at a constant speed of $\droneVelocity=6$\,m/s between locations, and takes $\timeAtLocation=60$\,s at each location to descend, deliver the package, and ascend. Batteries are given an energy density of $\energyDensity=650$\,kJ/kg and cost $\energyCostDrone=0.1$\,\$/kJ, based on the energy density and price of lithium polymer batteries in United States dollar values at the time this paper was written.

We use the \emph{percent improvement} $p$ to compare how a value changes from $x$ to $x'$ after adjusting a parameter or algorithm. We express the percent improvement as
\begin{equation}
p = \frac{x-x'}{x'}.
\end{equation}
We either minimize the cost or the overall delivery time, so when $p>0$, $x'$ is an improvement over $x$.

The rest of this section is organized as follows: \sref{ssec:hexacopter-energy-consumption-measurements} validates our energy consumption model by comparing it to measurements from a 3D Robotics ArduCopter Hexa-B hexacopter. \sref{ssec:mip-sa-performance} compares the results and runtimes of the MILP and SA DDP implementations for instances with 6-8 locations, and measures the consistency and runtimes of the SA DDPs for instances of up to 500 locations. \sref{ssec:limiting-budget-delivery-time} examines how the cost and overall delivery time found by the MC-DDP and MT-DDP respectively change with the delivery time limit, budget, area size, and number of locations. \sref{ssec:limiting-budget-delivery-time} also examines how the DDPs adjust costs for different delivery time limits and budgets. \sref{ssec:drone-reuse} investigates the importance of reusing drones, while \sref{ssec:battery-optimization} looks into the advantages of optimizing battery weight.

\subsection{Hexacopter Energy Consumption Measurements}
\label{ssec:hexacopter-energy-consumption-measurements}

In this section, we conduct two experiments to validate our energy consumption model. In the first experiment we verify that the power consumed during hover, level flight, and changing altitude is, on average, approximately equal. The second experiment verifies that the relationship between power consumption and battery and payload weight is approximately linear, and provides us with realistic values for our linear model.

Experimental results were generated with a 3D Robotics ArduCopter Hexa-B hexacopter. It has six 850 RPM/V motors with 10x4.6 APC slow fly propellers, uses an Arduino flight controller with firmware version 3.0.1 of the APM 2.5 autopilot, and has a u-blox LEA-6 GPS with a ceramic 4\,mm GPS patch antenna. Unless mentioned otherwise, the hexacopter is powered by a 4-cell 14.8\,V lithium polymer battery. Power consumption is calculated using current and voltage readings taken from the hexacopter's power management system, which were recorded by the ArduPilot mission planning software.

\begin{figure}[t!]
\centering
\subfloat[The hexacopter path during the straight flight manoeuvre.  Map data and image \copyright 2015 Google.\label{fig:straight-flight-path}]{\raisebox{0.25\height}{
		\includegraphics[width=\doubleFigureScaleValue]{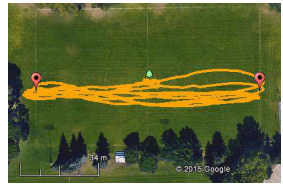}}
}
\hspace{20px}
\subfloat[Hexacopter altitude while performing altitude changes.\label{fig:alt-change}]{ \includegraphics[width=\doubleFigureScaleValue]{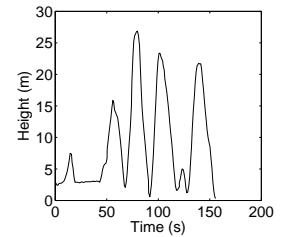}
}
		\caption{Hexacopter manoeuvre information.}
		\label{fig:manoeuvre-info}
\end{figure}

Power consumption was measured during level flight, hover, and altitude changes. The ArduPilot software under default settings was used to control the hexacopter and log the measurements. During level flight, the drone flew six times between two waypoints positioned about 94\,m apart, with an average speed of about 6\,m/s; the path can be seen in \fref{fig:straight-flight-path}. During hover, the hexacopter remained as stationary as it could over a single location over a period of 40\,s. During altitude changes, the hexacopter's altitude changed as in \fref{fig:alt-change}, with minimal changes in latitude and longitude.
\begin{table}[t!]
\caption{Average voltage, current, and power consumed by the hexacopter while performing different manoeuvres.}
\label{table:manoeuvre-power}
\footnotesize                                      
\centering                                                              
\begin{tabular}{|c|c c c|}                                        
\hline
Manoeuvre & Voltage (V) & Current (A) & Power (W) \\
\hline
Hover & 15.6 & 19.9 & 311  \\
Straight Flight & 15.3 & 19.8 & 304 \\
Altitude Change & 14.7 & 20.3 & 297 \\
\hline
\end{tabular}                                                                                                        
\end{table}

\tref{table:manoeuvre-power} contains the average voltage, current, and power consumption during hover, level flight, and altitude adjustment. It shows that, on average, the power consumed when performing the three manoeuvres is approximately equal. The current is slightly lower during level flight, possibly due to translational lift, where air moving horizontally along the rotors creates extra lift. The voltages in \tref{table:manoeuvre-power} are slightly different between tests, depending on the battery charge remaining. The altitude change tests had slightly lower voltages because they were conducted last.

The second set of experiments measure the power consumption of the hexacopter as it carries different weights, while hovering 3-5\,m above the ground at the same location. Measurements were repeated twice for each weight and battery type. The hexacopter was tested with both a 3-cell 11.1\,V lithium polymer battery and a 4-cell 14.8\,V lithium polymer battery.

\begin{figure}[t]
\centering
\includegraphics[width=\singleFigureScaleValue]{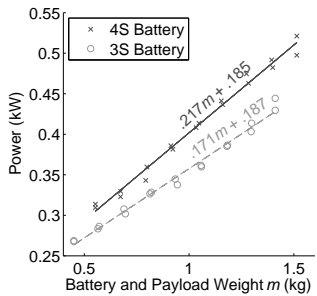}
\caption{Hexacopter energy consumption for various battery and payload weights $\batteryAndPayloadWeight$, with corresponding energy consumption model values.}
\label{fig:power-vs-weight}
\end{figure}

\fref{fig:power-vs-weight} confirms that the power consumption of the hexacopter varies approximately linearly with battery and payload weight, for both 3- and 4-cell batteries. The power consumed by the 4-cell battery grows faster with weight than with the 3-cell battery because we used 10x4.7 slow fly propellers, which perform less efficiently at the higher motor RPM from using the 4-cell batteries.

\subsection{Performance of MILP and SA Implementations}
\label{ssec:mip-sa-performance}

This section analyzes the performance of the MILP and SA implementations of the DDPs. It compares the results and runtimes of the SA implementations with the MILP implementations, and examines how the SA implementations perform at different cooling rate levels. The MILP implementations were run in CPLEX in a single thread, otherwise with default settings. Each scenario has its own number of locations $|\locationSetWithoutDepot{}|$, not including the depot, and area size $\areaSize$ in km$^2$.

\begin{table}[t!]        
\caption{Performance of MILP and SA implementations of the MT-DDP.}
\label{table:mip-sa-mt-dpp}                                               
\centering                                                         
\begin{tabular}{|c c|c c|c c c|}                                        
\hline                                 
\multicolumn{2}{|c|}{Scenario} & \multicolumn{2}{c|}{MILP} & \multicolumn{3}{c|}{SA} \\
$\areaSize$ & $|\locationSetWithoutDepot{}|$ & Result & Runtime & Mean & Std & Runtime \\
(km$^2$) &  & (min) & (s) & (min) & (min) & (s) \\
\hline
0.250 & 6 & 5.40 & 0.97 & 5.48 & 0.09 & 0.24 \\  
1.000 & 6 & 7.20 & 0.96 & 7.28 & 0.10 & 0.25 \\  
0.250 & 7 & 6.59 & 8.41 & 6.68 & 0.09 & 0.28 \\  
1.000 & 7 & 8.65 & 6.13 & 8.76 & 0.14 & 0.28 \\  
0.250 & 8 & 7.29 & 116.98 & 7.43 & 0.13 & 0.30 \\
1.000 & 8 & 9.63 & 86.87 & 9.83 & 0.22 & 0.31 \\
\hline
\end{tabular}                                                                                                        
\end{table}

\begin{table}[t!]        
\caption{Performance of MILP and SA implementations of the MC-DDP.}                                                
\label{table:mip-sa-mc-dpp}                                               
\centering                                                              
\begin{tabular}{|c c|c c|c c c|}                                        
\hline                                                                  
\multicolumn{2}{|c|}{Scenario} & \multicolumn{2}{c|}{MILP} & \multicolumn{3}{c|}{SA} \\
$\areaSize$ & $|\locationSetWithoutDepot{}|$ & Result & Runtime & Mean & Std & Runtime \\
(km$^2$) &  & (k\$) & (s) & (k\$) & (k\$) & (s) \\
\hline
0.250 & 6 & 0.99 & 9.47 & 1.00 & 0.00 & 0.19 \\   
1.000 & 6 & 1.04 & 4.84 & 1.04 & 0.00 & 0.19 \\   
0.250 & 7 & 1.04 & 75.38 & 1.04 & 0.00 & 0.22 \\  
1.000 & 7 & 1.10 & 24.84 & 1.11 & 0.00 & 0.22 \\  
0.250 & 8 & 1.04 & 1116.87 & 1.04 & 0.00 & 0.24 \\
1.000 & 8 & 1.22 & 205.52 & 1.26 & 0.04 & 0.24 \\ 
\hline
\end{tabular}                                                                                                        
\end{table}

Tables \ref{table:mip-sa-mt-dpp} and \ref{table:mip-sa-mc-dpp} show the average results and runtimes of the MILP and SA implementations of the MT-DDP and MC-DDP respectively. The SA algorithm parameters are set to $\initialTemperature=1$, $\finalTemperature=0.001$, $\temperatureAdjustment=0.9$, and $\temperatureRounds=1000$. MT-DDP scenarios have a budget $\budget$ of \$1,500, while MC-DDP scenarios have delivery time limits $\deliveryTimeLimit$ of 10 minutes.

The SA implementation consistently finds near-optimal solutions for problems with 8 or less locations. SA solutions to the MT-DDP are at most seconds away from optimal, while SA solutions to the MC-DDP are at most tens of dollars away, when optimal solutions are in minutes and thousands of dollars respectively. While optimality cannot be guaranteed when using SA, it appears to be an effective approach for solving small cases of the DDPs, providing similar results in under a second of runtime.

The runtime of the MILP implementation can grow exponentially with the number of locations, as can be observed from the runtimes in Tables \ref{table:mip-sa-mt-dpp} and \ref{table:mip-sa-mc-dpp}. In \tref{table:mip-sa-mc-dpp}, for example, the runtime increases from 75.38\,s to 1116.87\,s when solving an 8 location problem with an area of 0.25\,km$^2$ instead of a 7 location problem. The runtime of the SA implementation is consistently about 0.2-0.3\,s for each scenario.

In Tables \ref{table:mip-sa-mt-dpp} and \ref{table:mip-sa-mc-dpp}, increasing $\areaSize$ without changing $|\locationSetWithoutDepot{}|$ can actually cause the runtime of the MILP implementation to decrease. Locations are further apart, reducing the number of feasible routes and consequently the number of feasible solutions. Perhaps the solver improves the runtime by pruning solutions with routes that consume too much energy, carry too much weight, exceed the delivery time limit, or fail to meet some other criteria.

\newcommand{\columnSpacing}{1.05em}

\begin{table}[t]        
\caption{Simulated annealing results for the MT-DDP.}
\label{tab:sa-mt-dpp}                                               
\centering                                                              
\begin{tabular}{|c c | c @{\hspace{\columnSpacing}} c @{\hspace{\columnSpacing}} c | c @{\hspace{\columnSpacing}} c @{\hspace{\columnSpacing}} c|}                                        
\hline                                                                  
 &  & \multicolumn{3}{c|}{$\temperatureAdjustment=0.9$} &  \multicolumn{3}{c|}{$\temperatureAdjustment=0.99$}\\
$\areaSize$ & $|\locationSetWithoutDepot{}|$ & Mean & Std & Runtime & Mean & Std & Runtime \\
(km$^2$) &  & (min) & (min) & (s) & (min) & (min) & (s) \\
\hline                           
0.250 & 125 & 13.28 & 0.13 & 4.00 & 12.19 & 0.08 & 40.79 \\    
1.000 & 125 & 17.49 & 0.31 & 4.12 & 15.62 & 0.14 & 41.88 \\    
0.250 & 500 & 75.49 & 0.67 & 18.20 & 69.68 & 0.25 & 178.17 \\  
1.000 & 500 & 117.82 & 3.75 & 18.45 & 104.75 & 0.47 & 181.41 \\
\hline
\end{tabular}
\end{table}

\begin{table}[t]        
\caption{Simulated annealing results for the MC-DDP.}
\label{tab:sa-mc-dpp}                                               
\centering                                                              
\begin{tabular}{|c c | c c c |c c c|}                                        
\hline                                                                  
 &  & \multicolumn{3}{c|}{$\temperatureAdjustment=0.9$} & \multicolumn{3}{c|}{$\temperatureAdjustment=0.99$} \\
$\areaSize$ & $|\locationSetWithoutDepot{}|$ & Mean & Std & Runtime & Mean & Std & Runtime \\
(km$^2$) &  & (k\$) & (k\$) & (s) & (k\$) & (k\$) & (s) \\
\hline                           
0.250 & 125 & 13.99 & 0.23 & 5.82 & 13.52 & 0.21 & 57.91 \\  
1.000 & 125 & 17.03 & 0.29 & 6.04 & 16.21 & 0.24 & 58.22 \\  
0.250 & 500 & 60.40 & 1.02 & 49.31 & 54.55 & 0.43 & 461.27 \\
1.000 & 500 & 75.13 & 1.38 & 57.11 & 65.57 & 0.55 & 497.65 \\
\hline
\end{tabular}                                                                                                        
\end{table}

Tables \ref{tab:sa-mt-dpp} and \ref{tab:sa-mc-dpp} show the results of the SA implementation at different cooling rate levels $\temperatureAdjustment$ for scenarios with over 125 locations. As before, $\initialTemperature=1$, $\finalTemperature=0.001$, and $\temperatureRounds=1000$. The cooling rate level $\temperatureAdjustment$, however, is adjusted this time to see how much the results improve with additional runtime. MT-DDP scenarios have a budget $\budget$ of \$10,000, while MC-DDP scenarios have delivery time limits $\deliveryTimeLimit$ of 10 minutes.

The SA implementation results are consistent for large scenarios. In \tref{tab:sa-mt-dpp} the standard deviations are generally in seconds when the mean values are several minutes long. The standard deviation is always under \$1,500 in \tref{tab:sa-mc-dpp}, and usually just a few hundred dollars, where the means are often tens of thousands of dollars. The standard deviations are greatest when $|\locationSetWithoutDepot{}|$ is high and $\temperatureAdjustment$ is low. In this case the solution vector $\currentSolution$ is long, increasing the number of neighboring solutions that can be found after an adjustment, and $\temperatureAdjustment$ is low so the temperature decreases rapidly, giving fewer opportunities to make adjustments.

One possible reason for the small standard deviation is that the majority of costs appear to come from purchasing drones, which are considerably more expensive than energy. It seems like the SA DDPs consistently find similar numbers of drones, and differ more with regards to energy consumption. The variation in energy consumption is dampened by the relatively low energy costs, leading to the low standard deviations.

Increasing the cooling rate level from 0.9 to 0.99 resulted in percent improvements of up to 15\% for the means in the scenarios tested, which came at the cost of increasing runtime by about a factor of 10. Even though increasing the cooling rate level generally improves results, the improvement is not significant compared to the increase in runtime. While not shown in this paper, increasing $\temperatureAdjustment$ from 0.99 to 0.999 resulted in percent improvements of up to 5\%, again with the runtime increasing by about a factor of 10. We therefore set $\temperatureAdjustment=0.99$ for the rest of the section, as further increasing $\temperatureAdjustment$ would result in minor improvements compared to the increase in runtime.

We would like to mention that heuristics or exact algorithms that take advantage of characteristics unique to VRPs may find lower costs and times than the SA approach in this paper. An approach that takes advantage of geographical information to reduce the likelihood of connecting far apart locations would likely find better results as fewer infeasible routes would be tested. The number of times that performing operations such as swap, relocate, or 2-opt improve the solution could be tracked; operations which tend to improve the solution could be chosen more frequently than operations that do not.

While such approaches may lower cost, we do not believe that the trends we observed will change significantly with a different method. The trends can be explained by adjustments made to parameters such as area size and the number of customer locations in an area. They do not appear to be a result of behavior specific to SA.

\subsection{Analysis of DDP Parameters}
\label{ssec:limiting-budget-delivery-time}

This section examines how adjusting the budget $\budget$ in dollars, time limit $\deliveryTimeLimit$ in minutes, number of locations $|\locationSetWithoutDepot{}|$, and area size $\areaSize$ in km$^2$ affects the results obtained by the SA implementations of the DDPs. We show that the overall delivery time $\overallDeliveryTime$ found by the MT-DDP is inversely proportional to $\budget$, and that the cost $\totalCost$ found by the MC-DDP is inversely proportional to $\deliveryTimeLimit$. We also examine the sensitivity of the DDPs to changes in $\budget$, $\deliveryTimeLimit$, $|\locationSetWithoutDepot{}|$, and $\areaSize$. This information may be used to find the point where guaranteeing a better delivery time limit is not worth the additional cost, or alternatively the point at which reducing the budget is not worth the increase in overall delivery time.

\begin{figure}[t!]
\centering
\includegraphics[width=\singleFigureScaleValue]{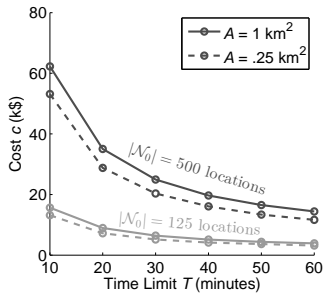}
\caption{Costs found by the SA MC-DDP algorithm for various delivery time limits, area sizes, and numbers of nodes.}
\label{fig:cost-vs-time}
\end{figure}

\begin{figure}[t!]
\centering
\includegraphics[width=\singleFigureScaleValue]{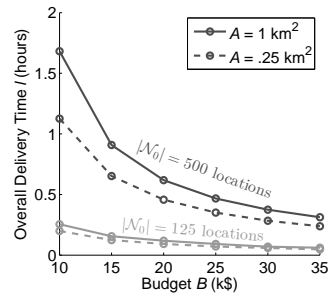}
\caption{Overall delivery times found by the SA MT-DDP algorithm for various budgets, area sizes, and numbers of nodes.}
\label{fig:time-vs-budget}
\end{figure}

\begin{figure}[t!]
\centering
\subfloat[MC-DDP cost breakdown\label{fig:mc-cost-breakdown}]{\includegraphics[width=\doubleFigureScaleValue]{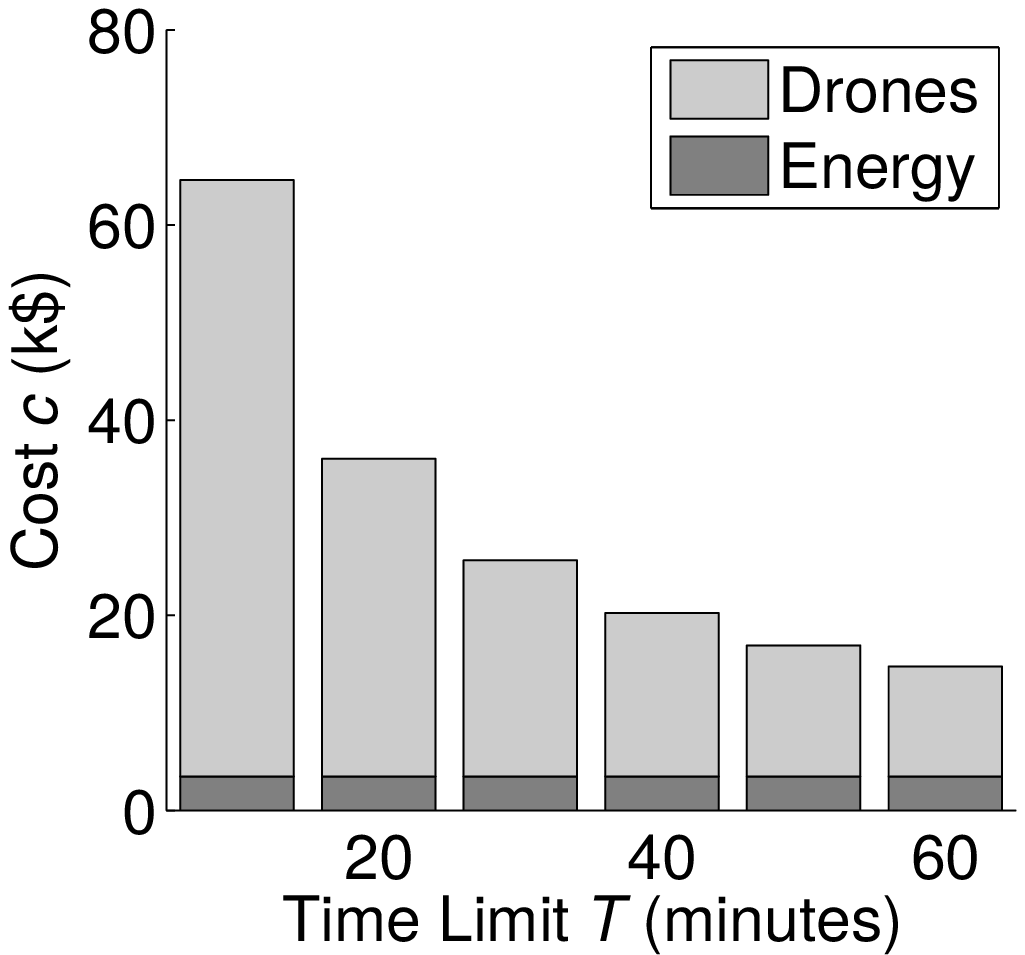}}
\subfloat[MT-DDP cost breakdown\label{fig:mt-cost-breakdown}]{\includegraphics[width=\doubleFigureScaleValue]{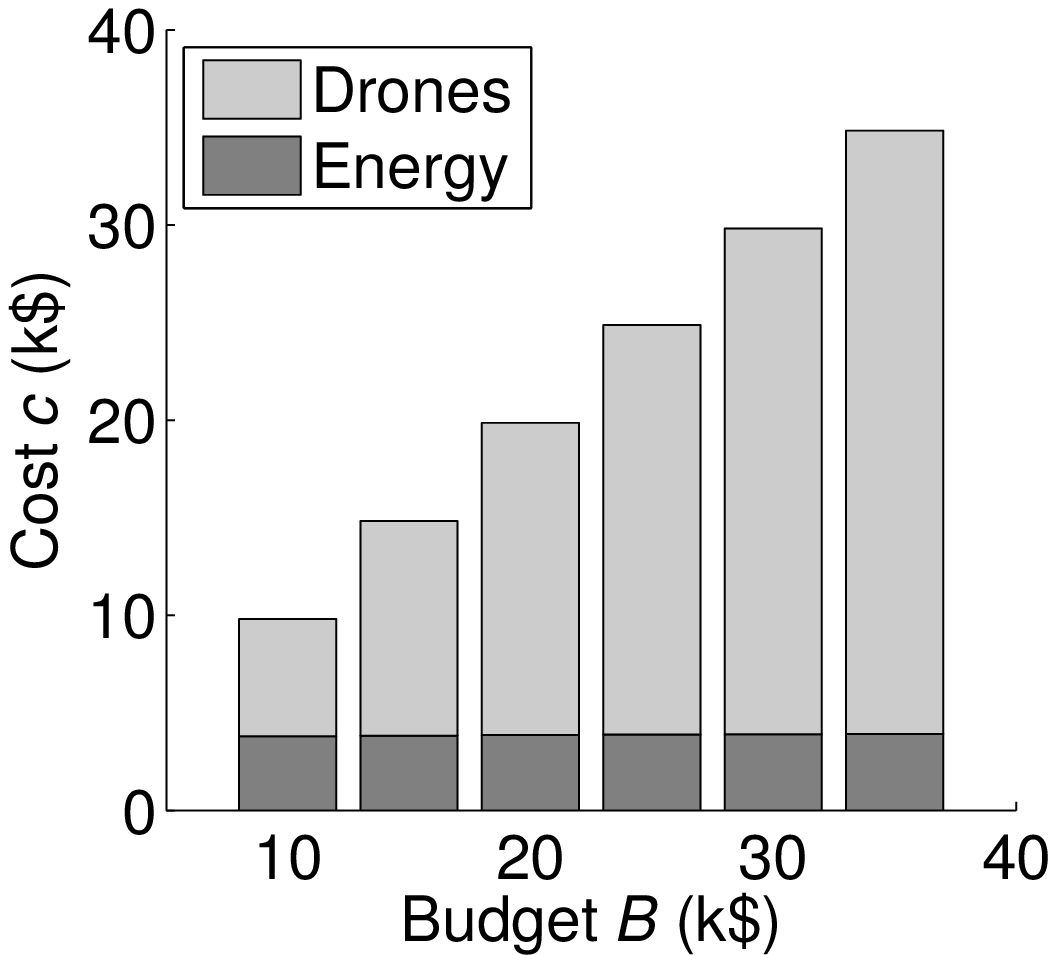}}
\caption{Costs found by the SA algorithms for various delivery time limits and budgets, with 500 locations in a 1\,km$^2$ area.}
\label{fig:cost-breakdown}
\end{figure}

\fref{fig:cost-vs-time} shows how adjustments to $\deliveryTimeLimit$, $|\locationSetWithoutDepot{}|$, and $\areaSize$ affect $\totalCost$, while \fref{fig:time-vs-budget} examines how adjustments to $\budget$, $|\locationSetWithoutDepot{}|$, and $\areaSize$ affect $\overallDeliveryTime$. To help explain the changes to $\overallDeliveryTime$ and $\totalCost$, Figs. \ref{fig:mc-cost-breakdown} and \ref{fig:mt-cost-breakdown} provide stacked bar graphs depicting the split in $\totalCost$ between drones and energy for scenarios with 500 locations in a 1\,km$^2$ area.

\fref{fig:cost-vs-time} shows an inverse exponential relationship between $c$ and $T$. Changes to $\totalCost$ from adjusting $\deliveryTimeLimit$ are larger when $\deliveryTimeLimit$ is low; from \fref{fig:mc-cost-breakdown} we can see that this change comes mostly from adding or removing drones. When $\deliveryTimeLimit$ is low, most drones fly short routes and are seldom reused; as $\deliveryTimeLimit$ increases, drone reuse and route length increase, while the number of routes decreases. This combination of factors reduces the number of purchased drones. After a certain point, increasing $\deliveryTimeLimit$ does not alter the routes as much, meaning costs are mainly reduced by reusing drones. The greater reductions in cost for low delivery times likely come from not only being able to reuse drones, but also from increasing the length and consequently reducing the number of routes flown.

\fref{fig:time-vs-budget} shows an inverse exponential relationship between $l$ and $B$. Changes to $\overallDeliveryTime$ from adjusting $\budget$ are greater when $\budget$ is low. From \fref{fig:mt-cost-breakdown} we can see that, as with the MC-DDP, this change comes mostly from adding or removing drones. When $\budget$ is low, only a small number of drones can be purchased to complete deliveries. Adding an additional drone reduces the delivery time considerably; for example, suppose a single drone is making all of the deliveries, then a second drone is added. If their burden is shared as evenly as possible, $\overallDeliveryTime$ can potentially be halved. Now suppose a large number of drones are available, and another one is added. The routes of each drone could perhaps be shortened slightly to give the additional drone a new route, but this will likely not have anywhere near the impact on $\overallDeliveryTime$ compared to adding a second drone when only a single drone is making deliveries.

Figs. \ref{fig:cost-vs-time} and \ref{fig:time-vs-budget} show that when $|\locationSetWithoutDepot{}|$ is high, the MC-DDP becomes more sensitive to changes in $\deliveryTimeLimit$, and the MT-DDP becomes more sensitive to changes in $\budget$. For a given increase in $\deliveryTimeLimit$ with the MC-DDP, more routes can be extended and additional drones can be reused when $|\locationSetWithoutDepot{}|$ is high. Similarly, increasing $|\locationSetWithoutDepot{}|$ when $\deliveryTimeLimit$ is low has a greater effect on cost than increasing $|\locationSetWithoutDepot{}|$ when $\deliveryTimeLimit$ is high. When $\deliveryTimeLimit$ is low, routes are short and drones are reused less, so increasing $|\locationSetWithoutDepot{}|$ requires purchasing more drones than if $\deliveryTimeLimit$ were high. For higher values of $|\locationSetWithoutDepot{}|$ with the MT-DDP, drones require more time to complete deliveries, as they must travel further and spend more time servicing locations. When $\budget$ is low, only a small number of drones are available and must therefore be reused more frequently, further increasing $\overallDeliveryTime$.

From Figs. \ref{fig:cost-vs-time} and \ref{fig:time-vs-budget} we see that the DDPs are less sensitive to changes in $\areaSize$ than to changes in $|\locationSetWithoutDepot{}|$. Increasing $\areaSize$ increases the distance traveled between locations, and consequently the time and energy required to complete deliveries. Changing $|\locationSetWithoutDepot{}|$ not only increases the travel distance, but also the time spent servicing locations, as well as the number of deliveries that must be made. Changing $|\locationSetWithoutDepot{}|$ therefore has a greater effect on time and cost than $\areaSize$. 

Note that these results may not hold for much higher values of $\areaSize$ or $|\locationSetWithoutDepot{}|$. If $\areaSize$ or $|\locationSetWithoutDepot{}|$ are so high that it becomes impossible in some scenarios to avoid violating the capacity constraint, cost constraint, or timing constraint, then the average price given by the SA algorithm could greatly increase. If a constraint is violated, the SA algorithm makes that solution undesirable by increasing its cost.

\subsection{Effect of Drone Reuse}
\label{ssec:drone-reuse}

In the previous section we learned that the majority of costs come from purchasing drones. Unlike other VRPs that model fuel consumption, the DDPs permit drones to be reused to perform multiple trips. To demonstrate the importance of allowing drone reuse, we show that our SA implementation performs better when drone reuse is permitted.

\begin{table}[t!]             
\centering                    
\caption{MC-DDP costs when reuse is enabled and disabled for various time limits in a 0.25\,km$^2$ area with 500 nodes.}      
\label{table:MCDDP-reuse}    
\begin{tabular}{|c|c c|c|}
\hline
 Delivery Time Limit & \multicolumn{2}{c|}{Cost (k\$)} & Percent Improvement  \\ 
(minutes) & Reuse & No Reuse &  (\%)  \\
\hline
10.0 & 54.6 & 112.9 & 106.84 \\
20.0 & 29.2 & 112.7 & 286.38 \\
30.0 & 20.6 & 112.8 & 446.65 \\
40.0 & 16.2 & 112.8 & 597.77 \\
50.0 & 13.6 & 112.7 & 730.95 \\
60.0 & 11.7 & 112.7 & 866.62 \\
\hline
\end{tabular}                 
\end{table}   

To look at the importance of drone reuse, we altered the SA implementation in \sref{ssec:route-costs} by removing the binary search in \aref{alg:drone-cost}, and setting the number of drones $\currentNumberOfDrones$ to the number of routes. \tref{table:MCDDP-reuse} compares the results of the MC-DDP with and without drone reuse. We did not include a similar table for the MT-DDP, as solutions could not be found for budgets of \$35,000 or less. 

\tref{table:MCDDP-reuse} shows substantial gains when allowing drone reuse. We can see that without reusing drones a limit of about \$113,000 is reached. By reusing drones, we can pass this limit without exceeding the delivery time constraints.

Note that we are ignoring the cost of maintaining drones. Flying fewer drones will likely increase maintenance costs due to increasing the total time flown per drone. New drones are generally expensive, and worn components can presumably be replaced quickly and relatively inexpensively, so reducing the number of drones while increasing the costs of labor and replacement parts is likely the less expensive option.

\subsection{Effect of Optimizing Battery Weight}
\label{ssec:battery-optimization}

While papers on minimizing fuel consumption do not reuse drones, MTVRPs do. We are not aware of any MTVRPs that try to optimize battery weight and therefore the energy given to each drone. Instead, they may limit each drone's travel distance or travel time to a predefined value. In this section, we apply our energy model to an MTVRP, and fix the battery weight for every route to the same constant value $\fixedBatteryWeight$ in kg. Similar to limiting distance or time, the amount of energy that can be consumed is limited. We show that optimizing battery weight for each route may provide better costs and delivery times compared to optimizing a fixed battery weight that is identical for each route. We also demonstrate that even if all routes share an identical battery weight, optimizing it can still be important.

\begin{figure}[t]
\centering
\includegraphics[width=\singleFigureScaleValue]{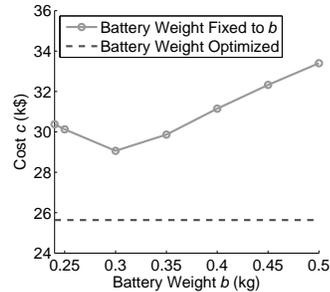}
\caption{Costs found by the MC-DDP when every drone's battery weighs $\fixedBatteryWeight$\,kg, with a time limit of 30 minutes on a 1\,km$^2$ field with 500 nodes.}
\label{fig:MCDDP-Weights}
\end{figure}

\begin{figure}[t]
\centering
\includegraphics[width=\singleFigureScaleValue]{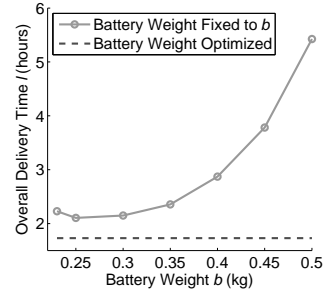}
\caption{Delivery times found by the MT-DDP when every drone's battery weighs $\fixedBatteryWeight$\,kg, with a budget of \$10,000 on a 1\,km$^2$ field with 500 nodes.}
\label{fig:MTDDP-Weights}
\end{figure}

To look at the effect of optimizing battery weight, we altered \aref{alg:energy-cost} to assume that every path has an identical battery weight of $\fixedBatteryWeight$\,kg. The energy in a battery is then $\energyDensity \fixedBatteryWeight$, where $\energyDensity$ is the energy density of the battery in kJ/kg. \fref{fig:MCDDP-Weights} shows the cost $\totalCost$ found by the MC-DDP when assuming an identical battery weight $\fixedBatteryWeight$ for each route. The dashed line represents the cost found by optimizing battery weights for each route with the unaltered version of \aref{alg:energy-cost}; it is flat because it does not consider $\fixedBatteryWeight$. \fref{fig:MTDDP-Weights} is similar, but for the overall delivery time $\overallDeliveryTime$ found by the MT-DDP.

Figs. \ref{fig:MCDDP-Weights} and \ref{fig:MTDDP-Weights} demonstrate that optimizing battery weights for each route to give drones exactly enough energy to complete their upcoming routes can provide substantial savings compared to optimizing a fixed battery weight that is identical for every route. In \fref{fig:MCDDP-Weights} the percent improvement when optimizing each route's battery weight amounts to 13\% compared to assuming an identical battery weight of 0.3\,kg for every route. In \fref{fig:MTDDP-Weights} the percent improvement of optimizing each route's battery weight amounts to about 22\% compared to assuming each route has a battery weight of 0.25\,kg when using the MT-DDP. Note that while optimizing battery weight for each route provides savings compared to optimizing a fixed battery weight, it can be less effective with a small area size or low number of nodes.

Even if the battery weight is identical for all routes, optimizing it can improve cost and delivery time, as can be seen in Figs. \ref{fig:MCDDP-Weights} and \ref{fig:MTDDP-Weights}. Setting $\fixedBatteryWeight$ to 0.3\,kg instead of 0.5\,kg provides a percent improvement of 15\% for the MC-DDP, and setting $\fixedBatteryWeight$ to 0.25\,kg instead of 0.45\,kg provides a percent improvement of 80\% for the MT-DDP.

\section{Conclusion}
\label{sec:conclusion}

In this paper, we derived an energy consumption model for multirotor drones, provided a linear approximation for it, and verified the approximation experimentally. With this approximation we proposed the drone delivery problems (DDPs) that minimize cost or delivery time while considering battery weight, payload weight, and drone reuse, and implemented them as mixed integer linear programs. To solve practical scenarios with hundreds of locations, we proposed a string-based simulated annealing algorithm for solving the DDPs.

Numerical results indicate that optimizing battery weight and reusing drones are important considerations for drone delivery. Optimizing battery weights resulted in percent improvements of over 10\% compared to solutions where each drone had an identical battery weight. Even when drones had an identical battery weight, optimizing that weight provided a percent improvement of 80\% for the MT-DDP in the scenario considered in \sref{ssec:battery-optimization}. In \sref{ssec:drone-reuse} reusing drones led to costs in the tens of thousands of dollars, while preventing drone reuse resulted in costs of approximately \$113,000.

When minimizing costs with the MC-DDP, an inverse exponential relationship between the cost and the delivery time limit was noticed. An inverse exponential relationship was also seen between the overall delivery time and the budget constraint with the MT-DDP. For both DDPs, the majority of costs are allocated to drones instead of energy. This information is useful for a practitioner. For example, reducing already low delivery times can require a relatively large investment in additional drones. Adding a single drone when the number of drones is low can greatly reduce delivery times, while adding a drone when the number of drones is high has a limited effect in comparison.

One weakness of the SA algorithm is that it does not take advantage of characteristics inherent to a VRP. For instance, it does not take advantage of geographical information to reduce the likelihood of trying infeasible routes with two locations at opposite ends of the area of interest. The cost function of the SA algorithm accounts for reusing drones, as well as the effect of battery and payload weight on energy consumption, so it could be applied to a heuristic better suited to VRPs. The cost function itself is also an approximation, as it uses list scheduling to assign drones to routes. Perhaps a multi-layer approach that optimizes routes and their assignment to drones would perform better.

Additional strategies could be implemented to further reduce costs and generalize the problem. Time windows could be added to locations to ensure that packages are delivered within a specific time. Multiple depots or recharging stations could be studied as a method of extending the flight distance of drones. While we treat cost and delivery time as separate optimization objectives, a multi-objective optimization problem could combine both. Maintenance costs could be included when determining drone reuse. The performance of the MILP and SA DDPs on a cloud service such as Amazon Web Services \cite{aws} could be investigated. The runtime or resources provided to the algorithms could be adjusted to compare performance gains from increasing the service cost.

The DDPs have practical applications in delivery and emergency response scenarios. They can plan routes and assign drones to them in scenarios with tight budgets or time limits. In a disaster response scenario, for example, the MT-DDP could make the best use of a limited budget to optimize the routes taken by a fleet of drones to deliver food or medicine to assigned locations as quickly as possible. Alternatively, a delivery company could determine the budget necessary to guarantee a minimum delivery time with the MC-DDP. The SA implementations can be used to quickly predict delivery times and plan budgets for drone-based delivery operations with a large number of locations.

\appendices

\bibliographystyle{IEEEtran}
\bibliography{hardware,technology,vrp,aerodynamics,optimizers,misc}

\newcommand{\biospace}{-0.5in}
\vspace{\biospace}
\begin{IEEEbiography}
[{\includegraphics[width=1in,height=1.25in,clip,keepaspectratio]{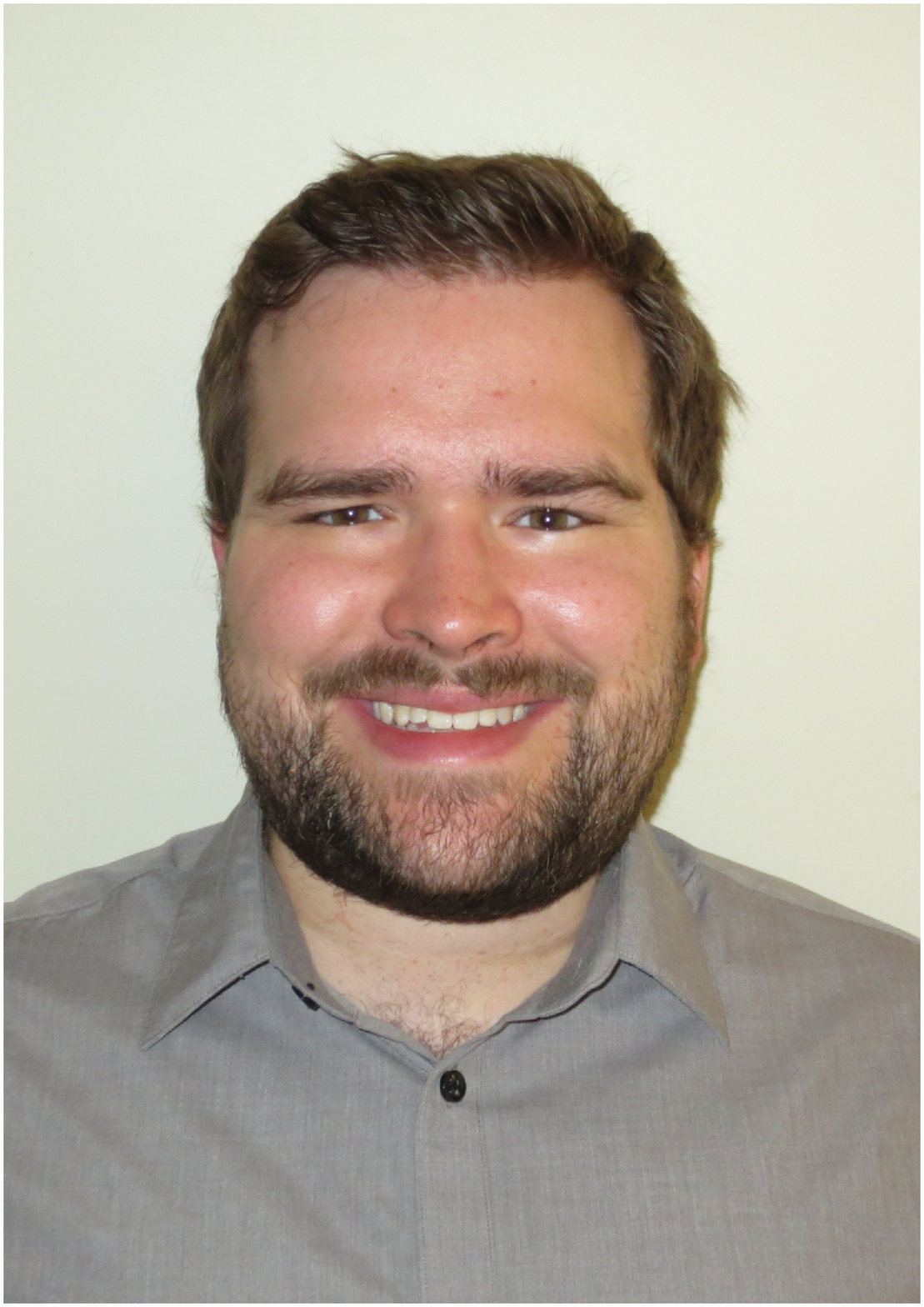}}]
{Kevin Dorling}
(S'14) received a B.S. degree in Computer Engineering with distinction from the University of Calgary, Canada, in 2010. He is currently pursuing a Ph.D. in Electrical Engineering at the University of Calgary.

In 2011 and 2013, he was an intern at Bell Labs, Alcatel-Lucent in Stuttgart, Germany, working on deployment techniques for wireless sensor networks. In 2009, he was an intern at CDL Systems in Calgary, Canada, developing control station software for unmanned vehicles. His research interests include wireless sensor networks, operations research, and communication systems.
\end{IEEEbiography}
\vspace{\biospace}
\begin{IEEEbiography}
[{\includegraphics[width=1in,height=1.25in,clip,keepaspectratio]{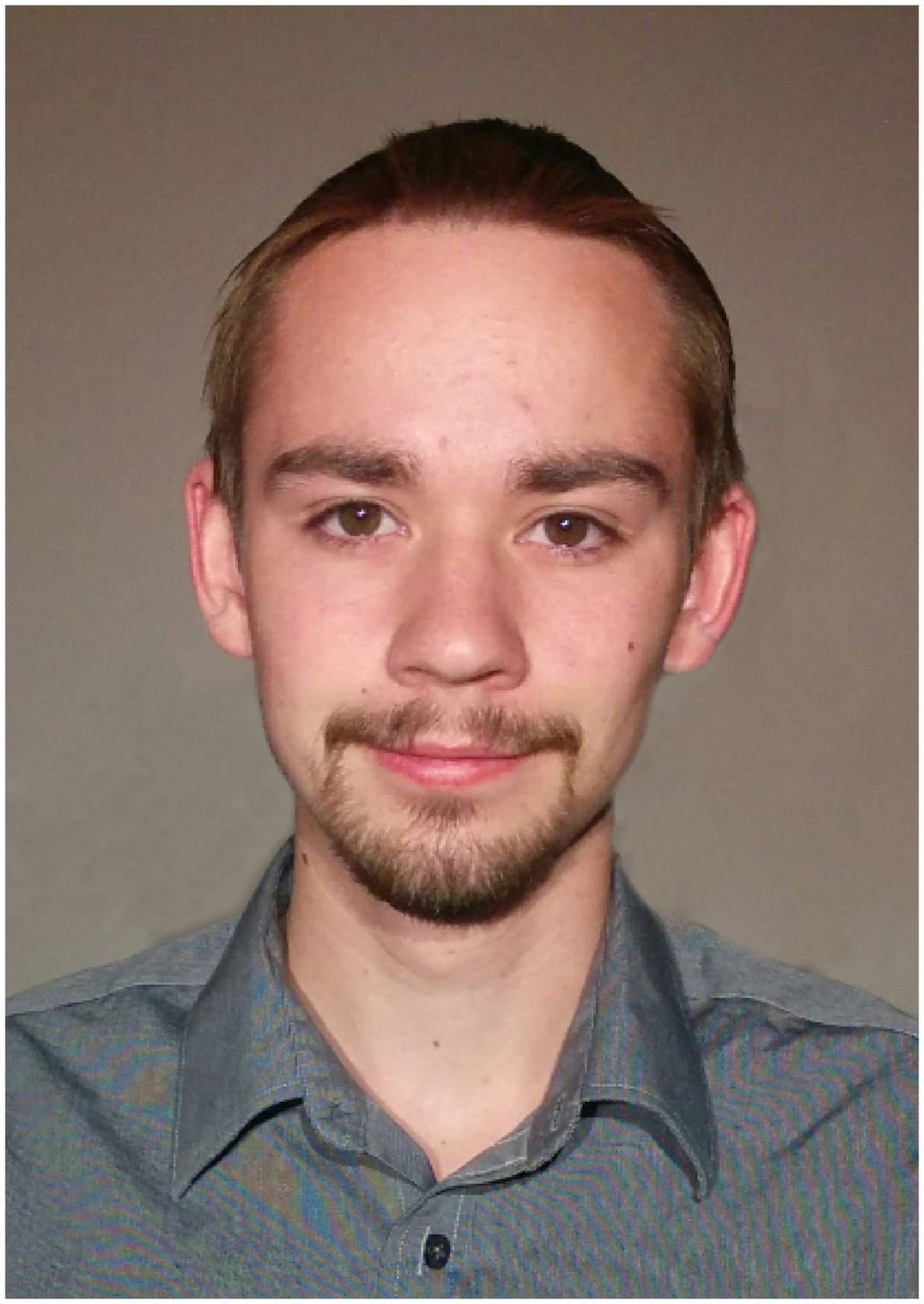}}]
{Jordan Heinrichs}
received a B.S. degree in Electrical Engineering, minor in Computer Engineering with distinction from the University of Calgary, Canada, in 2016.

In 2013 he tested multicopters as a summer research student. From 2014 to 2015 he interned at Lockheed Martin CDL Systems developing software for miniature unmanned aerial vehicles. He is currently employed at Hitachi ID Systems developing software for identity management.
\end{IEEEbiography}
\vspace{\biospace}
\begin{IEEEbiography}
[{\includegraphics[width=1in,height=1.25in,clip,keepaspectratio]{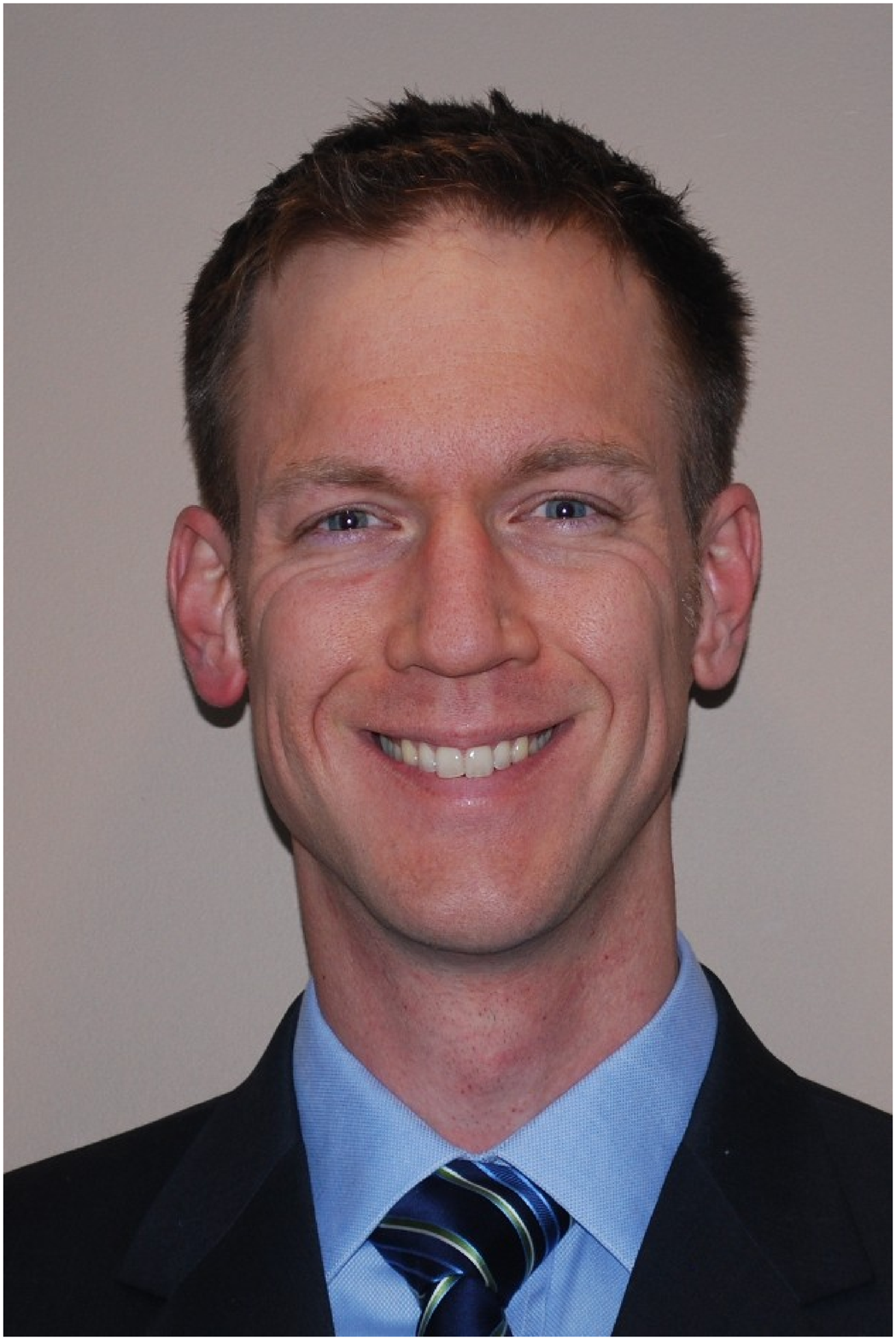}}]
{Geoffrey Messier}
(S'91 - M'98) received his B.S. in Electrical
Engineering and B.S. in Computer Science degrees from the University of
Saskatchewan, Canada with great distinction in 1996.  He received his
M.Sc. in Electrical Engineering from the University of Calgary, Canada
in 1998 and his Ph.D. degree in Electrical and Computer Engineering from
the University of Alberta, Canada in 2004.

From 1998 to 2004, he was employed in the Nortel Networks CDMA Base
Station Hardware Systems Design group in Calgary, Canada.  At Nortel
Networks, he was responsible for radio channel propagation measurements
and simulating the physical layer performance of high speed CDMA and
multiple antenna wireless systems.  Currently, Dr. Messier is a
Professor in the University of Calgary Department of
Electrical and Computer Engineering.  His research interests include
data networks, physical layer communications and communications
channel propagation measurements.
\end{IEEEbiography}
\vspace{\biospace}
\begin{IEEEbiography}
[{\includegraphics[width=1in,height=1.25in,clip,keepaspectratio]{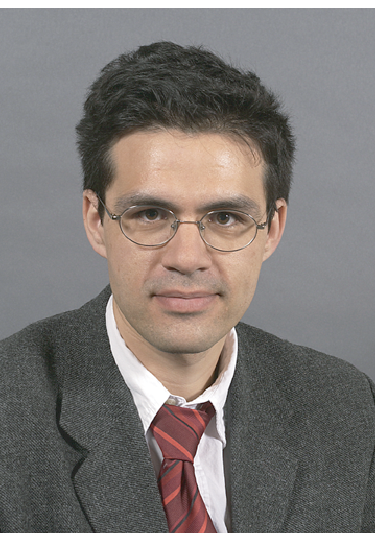}}]
{Sebastian Magierowski} received his Ph.D. degree in Electrical Engineering from the University of Toronto in 2004.  From 2004 to 2012 he served as Assistant/Associate Professor in the Department of Electrical and Computer Engineering at the University of Calgary after which he joined the faculty of the Department of Electrical Engineering and Computer Science in the Lassonde School of Engineering at York University, Toronto, Canada.  

As part of his industrial experience (Nortel Networks, PMC-Sierra, Protolinx Corp.) Dr. Magierowski has worked on CMOS device modeling, high-speed mixed-signal IC design, and data networks.  His research interests include analog/digital CMOS circuit design, communication systems and biomedical instrumentation.
\end{IEEEbiography}

\end{document}